\documentclass[12pt]{article}
\usepackage{amssymb,amsmath,amsthm, amsfonts}
\usepackage{graphicx}
\usepackage{epsfig}
\def\dref#1{(\ref{#1})}

\theoremstyle{plain}
\newtheorem{theorem}{Theorem}[section]
\newtheorem{lemma}{Lemma}[section]

\newtheorem{corollary}{Corollary}[section]

\theoremstyle{definition}

\numberwithin{equation}{section}

\begin{document}

\title{\large\bf Asymptotic profile  of  a two-dimensional chemotaxis--Navier--Stokes system
 with singular sensitivity and logistic source}

\author
{\rm Peter Y.~H.~Pang\footnote{Email: {\tt matpyh@nus.edu.sg}}\\
\it\small Department of Mathematics, National University of Singapore\\
\it\small 10 Lower Kent Ridge Road, Republic of Singapore 119076
\\[3mm]
\rm Yifu Wang\footnote{Corresponding author. Email: {\tt wangyifu@bit.edu.cn}}\\
\it\small School of Mathematics and Statistics, Beijing Institute of Technology\\
\it\small Beijing 100081, People's Republic of China
\\[3mm]
\rm Jingxue Yin\footnote{Email: {\tt yjx@scnu.edu.cn}}\\
\it\small School of Mathematics, South China Normal University\\
\it\small Guangzhou, Guangdong 510631, People's Republic of China}

\date{}

\maketitle \vspace{0.3cm}
\noindent
\begin{abstract}
 The chemotaxis--Navier--Stokes system
\begin{equation*}\label{0.1}
 \left\{\begin{array}{ll}
  n_t+u\cdot \nabla n=\triangle n-\chi\nabla\cdotp \left(\displaystyle\frac n {c}\nabla c\right)+n(r-\mu n),\\
 c_t+u\cdot \nabla c=\triangle c-nc,\\
 u_t+ (u\cdot \nabla) u=\Delta u+\nabla P+n\nabla\phi,\\
 \nabla\cdot u=0,
  \end{array}\right.
\end{equation*}
is considered in a bounded smooth domain $\Omega \subset \mathbb{R}^2$, where $\phi\in W^{1,\infty}(\Omega)$, $\chi>0$, $r\in \mathbb{R}$ and $\mu> 0$ are given  parameters.
It is shown that  there exists  a  value
$\mu_*(\Omega,\chi, r)\geq 0$ such that whenever $ \mu>\mu_*(\Omega,\chi, r)$,
the  global-in-time classical solution
to the system is  uniformly bounded with respect to $x\in \Omega$. Moreover, for the case $r>0$,
$(n,c,\frac {|\nabla c|}c,u)$ converges to $(\frac r \mu,0,0,0)$ in $L^\infty(\Omega)\times L^\infty(\Omega)\times L^p(\Omega)\times L^\infty(\Omega)$ for any $p>1$ exponentially as $t\rightarrow \infty$, while in the case $r=0$, $(n,c,\frac {|\nabla c|}c,u)$  converges to $(0,0,0,0)$ in $(L^\infty(\Omega))^4$
  algebraically.
 To the best of our knowledge, these results provide the first precise information on the asymptotic profile of solutions %to \eqref{0.1}  in two dimensions.

\end{abstract}

\vspace{0.1cm}
%\begin{spacing}{2.0}
\noindent {\bf\em Keywords:}~Chemotaxis, Navier--Stokes, singular sensitivity, asymptotic profile.

\noindent {\bf\em AMS Classification \rm (2010)}:~35K55, 35Q30, 35B40, 35Q92, 76D05, 92C17.
\section{Introduction}

In their seminal work (\cite{Keller2,Keller}), using cross-diffusive parabolic PDE systems, Keller and Segel studied the phenomenon of the emergence of spatial structures in biological systems through chemically induced processes. In particular, they looked at situations where componenets of the biological system were able to actively secrete a chemoattractive signal which then directed the movements of components in the system; or where, instead of having the ability to produce such signals, components of the system simply followed and consumed a chemical nutrient produced externally. A prototypical example of the former is the dictyostelium discoideum colony, while that of the latter is an E. coli population. The movement thus induced by either a chemical signal or nutrient is called chemotaxis, and the corresponding  model of the latter type is sometimes called a chemotaxis--consumption model.
Often, such chemotactic movements take place in a fluid environment, and experimental findings and analytical studies have revealed  the remarkable effects of  %liquid environment of populations %
chemotaxis--fluid interaction
on the overall behavior of the respective chemotaxis systems, such as the prevention of blow-up and improvement of efficiency of mixing (\cite{Chertock, Kiselev2, Kiselev3, Lorz,Tuval}).

In this paper, we are concerned with  the
 chemotaxis--consumption system coupled with the incompressible Navier--Stokes equations
\begin{equation}\label{1.1}
 \left\{\begin{array}{ll}
  n_t+u\cdot \nabla n=\triangle n-\chi\nabla\cdotp \left(\displaystyle\frac n {c}\nabla c\right)+f(n),\,&
x\in \Omega, t>0,\\
 c_t+u\cdot \nabla c=\triangle c-nc, \,&
x\in \Omega, t>0,\\
 u_t+ (u\cdot \nabla) u=\Delta u+\nabla P+n\nabla\phi,\,&
x\in \Omega, t>0,\\
 \nabla\cdot u=0,\,&
x\in \Omega, t>0,\\
  \end{array}\right.
\end{equation}
describing the biological population density $n$,  the chemical signal concentration  $c$, the incompressible fluid velocity $u$ and the associated pressure  $P$ of the fluid flow
in the physical domain $\Omega\subset\mathbb{R}^N$. It is assumed  that  $n$ and $c$ diffuse randomly as well as are transported by the fluid, with a bouyancy effect on $n$ through the presence of a given gravitational potential $\phi$. Further, it is assumed that the chemotactic stimulus is perceived in accordance with the Weber--Fechner Law (\cite{Short,Zwang0,WinklerAIP}) which states
that subjective sensation is proportional to the logarithm of the stimulus intensity, in other words,
the population $n$ partially direct their movement toward
increasing concentrations of the chemical nutrient $c$ that they consume with the logarithmic sensitivity.  In addition, on the considered time scales of cell migration,
we allow for population growth to take place, through the term $f(n)=rn-\mu n^2$ with the effective growth rate $r\in \mathbb{R}$ and strength of the overcrowding effect $\mu>0$;
we note that $r=0$ is allowed and has indeed been argued for in certain models (\cite{Hillen,Kiselev2}).

The system \eqref{1.1} appears to generate interesting, non-trivial dynamics. However,
to the best of our knowledge, no analytical result is available yet which rigorously describes the qualitative behavior of such
solutions.
This may be due to the circumstance that \eqref{1.1} joins two subsystems which are far from being fully understood even when decoupled from each other. Indeed, %considered separately.
  %as a subsystem,
 \eqref{1.1} contains  the Navier--Stokes
equations which themselves do not admit a complete existence and regularity
theory (\cite{Wiegner}).
%and elaborate research since Leray��s pioneering work (cf. for example [30] for a survey, and also [10,24,28]). On the other hand, the first

At the same time, by setting $u\equiv 0$ in \eqref{1.1},
%As for the latter, it
we arrive at the following chemotaxis--consumption model
%without the population growth %when $f(n)=0$
\begin{equation}\label{1.2}
 \left\{\begin{array}{ll}
  n_t=\triangle n-\chi \nabla\cdotp (\frac {n}c\nabla c),\\
 c_t=\triangle c-nc,
  \end{array}\right.
\end{equation}
where population growth has been ignored, which was introduced by Keller and Segel (\cite{Keller}) to describe the collective behaviour of the bacteria E.\ coli set in one end of a
capillary tube featuring a gradient of nutrient concentration observed in
the celebrated experiment of Adler (\cite{Adler}). Later, this model was also  employed to describe the dynamical interactions between vascular endothelial cells and vascular
endothelial growth factor (VEGF) during the initiation of tumor angiogenesis (see \cite{Corrias, Levine}).
It has already been demonstrated that the logarithmic sensitivity featured in \eqref{1.2} renders a significant degree of complexity in the system; in particular, it plays an indispensable role in  generating wave-like solutions without any type of cell kinetics (\cite{Hillen,Keller2,Rosen,Schwetlick, Zwang0}), which
is  a prominent feature in the Fisher equation (\cite{Kolmogorov}).

In comparison with \eqref{1.2}, the related chemotaxis system %signal--production relative thereof %of \eqref{1.1}
\begin{equation}\label{1.3}
 \left\{\begin{array}{ll}
  n_t=\triangle n-\chi \nabla\cdotp (\frac {n}c\nabla c)+f(n),\\
 c_t=\triangle c-c+n,
  \end{array}\right.
\end{equation}
where the chemical signal $c$ is actively secreted by the bacteria rather than consumed (see \cite{Bellomo,Hillen}), has been more extensively studied.
It is observed that the chemical signal production mechanism in the $c-$equation inhibits the tendency of %  an evolution of
$c$ to take on small values, and thereby the singularity in the sensitivity function is mitigated. Accordingly, for such higher dimensional systems with reasonably smooth but
arbitrarily large data,  the global existence of bounded smooth solutions can be achieved.
Indeed, global existence and boundedness of classical solutions to \eqref{1.3} without source terms %in absence of  sources %with  $f(n)=0$
 is guaranteed if $\chi\in(0,\sqrt{\frac2N})$ (\cite{Fujie,Winklermmas}), or if
$N=2, \chi\in(0,\chi_0)$ with some $\chi_0>1.015$ (\cite{Lankeit}), while certain generalized solutions have been constructed
for general $\chi>0$ in the two-dimensional radially symmetric case (\cite{Stinner,Winklermmas}). Moreover, without any symmetry hypothesis,
 % the requirement of  radial symmetry,
 Winkler and Lankeit established the global solvability of generalized solutions for the cases $\chi<\infty, N=2$; $\chi<\sqrt{8}, N=3$; and $\chi<\frac N{N-2}, N\geq 4$  (\cite{LankeitNdea}).

% seem appropriate, and  a prototypical choice of $f$ is the standard logistic type $n(r-\mu n) $.
Furthermore, in accordance with known results for the classical Keller--Segel chemotaxis model (see \cite{Lankeit15, Winklercpde, WinklerJDE14} for example), the presence of the logistic source term $f(n)=n(r-\mu n)$ in \eqref{1.3}
can inhibit the tendency toward explosions of cells at least under some restrictions on certain parameters. Indeed, it is known that
\eqref{1.3} with $N = 2$ possesses
a global classical solution $(n,c)$ for any $r\in \mathbb{R}, \chi, \mu>0$, and $(n,c)$ is globally bounded if
$r>\frac{\chi^2} 4$ for $0 < \chi\leq 2$ or $r > \chi-1$ for $\chi>2$ (\cite{Zhao1}). Moreover, $(n,c)$ exponentially converges to $(\frac r{\mu},\frac r{\mu})$ in $L^\infty(\Omega)$ %as $t\rightarrow \infty $
provided that $\mu>0$ is sufficiently large (\cite{Zheng}). As for the higher dimensional cases ($N\geq 2$),  the global very weak solution of \eqref{1.3} with
 $f(n)=r n-\mu n^k$ is constructed when $k,\chi$ and $r$ fulfill a certain condition. In addition, when $N=2$ or $3$, this solution is global bounded provided  $\frac r\mu$ and the initial data
$\|n_0\|_{L^2}, \|\nabla c_0\|_{L^4}$ are suitably small (\cite{Zhao3}).

In contrast to \eqref{1.3},
 system \eqref{1.2} is more challenging
due to the combination of the consumption of $c$ with
the singular chemotaxis sensitivity of $n$.  Intuitively, the absorption mechanism in the $c-$equation of \eqref{1.2},
which induces the preference for small values of $c$, considerably intensifies the destabilizing potential of singular sensitivity
in the $n-$equation.
Up to now, it seems that only limited results on global  classical solvability in the spatial two-dimensional case are available.
 In fact, only recently have certain global generalized solutions to
\eqref{1.2} been constructed for general initial data in \cite{Lankeit2, Winkler3M, WinklerJDE}, whereas with respect to global classical solvability, %of the associated boundary value problems,
it has only been shown %is restricted to
for some small initial data (see \cite{Zwang,Winklerp}). In particular, Winkler (\cite{Winklerp}) showed that the global classical solutions to \eqref{1.2} in bounded convex two-dimensional domains exist and  converge to the homogeneous steady state under an essentially explicit smallness condition on $n_0$ in $L\log L(\Omega)$ and $\nabla\ln c_0 $ in $L^2(\Omega)$.

We would, however, like to note that numerous variants of \eqref{1.2}, such as those involving nonlinear diffusion, logistic-type cell kinetics and saturating signal production
 %considered recently. and  �� various biological mechanisms involving the porous medium-type nonlinear diffusion, logistic type cell kinetics and saturating signal production
%have been considered recently. In these situations, the corresponding variants of \eqref{1.2} possess global classical solutions for rather
%general initial data
(\cite{Ding,Jia,Lankeitnon,Lankeitjde,Lankeitv,Liu,Viglialoro,Zhao2}), have been studied. For example, the authors of \cite{Zhao2} proved that the particular version of \eqref{1.2} by adding  %$S(u)=u, % ^\beta ~(0<\beta<1),
 $f(n)=rn-\mu n^{k}~ (r>0, \mu>0, k>1)$ into the $n$-equation admits a  global classical solution $(n,c)$ in the bounded domain $\Omega\subset \mathbb{ R}^N$ if $k>1+\frac N2$,
and in the two-dimensional setting,  $(n,c, \frac {|\nabla c|}c)\rightarrow  ((\frac r{\mu})^{\frac 1{k-1}},0,0)$ for sufficiently large $\mu$.
In particular, it is shown in the recent paper \cite{Lankeitnon} that \eqref{1.2} with  logistic source $f(n)=r n-\mu n^2~ (r\in \mathbb{R}, \mu>0)$  possesses a unique global classical solution %$(u,v)$
if $0<\chi<\sqrt\frac{2}N$, $\mu>\frac{N-2}{2N}$, and a globally bounded solution only in \bf{one}   \rm  dimension for any $\chi>0, \mu>0$.
Also, the author of \cite{w} showed that if $\mu>\mu_0$ with some $\mu_0=\mu_0(\Omega,\chi)>0$ %is linearly bounded from below by $r$,
then the corresponding classical solution is globally bounded, and
 $(n,c, \frac {|\nabla c|}c)\rightarrow  (\frac {r_+}{\mu},\lambda,0)$  with $\lambda\in [0,\frac 1 {|\Omega|}\int_\Omega c_0)$ in
$(L^\infty(\Omega))^3$ as $t\rightarrow\infty$. Of course, this leaves  open the possibility  of  blow-up of solutions when $\mu$ is positive  but small. Anyhow, it have been shown in  \cite{WinklerDCDSB2017} that when $\mu>0$ is suitably small, the strongly destablizating action of chemotactic cross--diffusion may lead to the occurrence of solutions which attain possibly finite but arbitrarily large values.

Coming back to our chemotaxis--consumption--fluid model \eqref{1.1}, as we have already pointed out, very little seems to be known regarding
%In comparison to the latter,  %chemotaxis--fluids system,  %o the best of our knowledge,
%very few seems to be known  %the  study of  chemotaxis--fluid systems with  logarithmic sensitivity seems  very few so far
%
the qualitative behavior of solutions
%in the framework of  chemotaxis--fluid systems with  logarithmic sensitivity, inter alia \eqref{1.1}
(\cite{Black,BlackJEQ,BlackMMA}).
In fact, we are aware of one result only which is concerned with the asymptotic behavior and
eventual regularity of solutions to  the Stokes--variant of \eqref{1.1}. Namely, it is shown in \cite{Black} that for small initial mass $ \int_\Omega n_0$,
the corresponding system upon neglection of $u\cdot \nabla u$ and $f(n)$ in \eqref{1.1} possesses at least one global generalized solutions, which will become
smooth after some waiting time and stabilize toward the steady state $ (\frac 1 {|\Omega|}\int_\Omega n_0,0,0)$ with respect to the topology of
$(L^\infty(\Omega))^3$.
 %obey certain asymptotic propertie
% ment In contrast to this, the  knowledge  on the chemotaxis--fluid system
Since the presence of the fluid interaction does not have any regularizing effect on the large time behavior, it is expected
that instead of the small restriction on the initial data, the quadratic degradation %in \eqref{1.1}
may have a substantial regularizing effect on the dynamic behavior of solutions to \eqref{1.1}. %the corresponding chemotaxis Navier--Stokes frameworks.

%%%%%%
%Accordingly, it is quite natural to ask whether the solutions of \eqref{1.1} are uniformly bounded with respect to time, and especially what their  large time behavior is when  the coefficient $\mu$ is suitably large as compared to $r$, or when $r=0,\mu>0$.

The goal  of the present work is to give the asymptotic profile in time %fairly complete description
of solutions to \eqref{1.1}
in the two-dimensional case. In order to state our main results, we shall impose on \eqref{1.1}
the boundary conditions
\begin{equation}\label{1.4}
\nabla n\cdot \nu=\nabla c\cdot \nu=0\quad\hbox{and}~~u=0 ~~ \hbox{for}~~ x\in\partial\Omega,
\end{equation}
and initial conditions
\begin{equation}\label{1.5}
n(x,0)=n_0(x),~ c(x,0)=c_0(x), ~  u(x,0)=u_0(x)~ ~\hbox{for}~~x\in \Omega.
\end{equation}
Throughout this paper, it is assumed that
\begin{equation}\label{1.6}
\left\{\begin {array}{l}
n_0\in C^0(\bar{\Omega}), ~n_0\geq 0~ \hbox{and}~ n_0\not\equiv 0 ~ \hbox{in}~ \Omega,\\
c_0\in W^{1,\infty}(\Omega), c_0>0 ~ \hbox{in}~ \bar{\Omega}~~\hbox{as well as} \\
u_0\in D(A^{\beta})\;\hbox{for all}\;\beta\in(\frac12,1)
\end{array}\right.
\end{equation}
with $A$ denoting the Stokes operator $A=-\mathcal{P}\Delta$ with domain $D(A):=W^{2,2}(\Omega; \mathbb{R}^2)\cap W^{1,2}_0(\Omega; \mathbb{R}^2)\cap L_\sigma^2(\Omega)$,
where $L_\sigma^2(\Omega):=\{\varphi\in L^2(\Omega; \mathbb{R}^2)|\nabla\cdot\varphi=0\}$ and $\mathcal{P}$ stands for the Helmholtz projection of $L^2(\Omega)$ onto $L_\sigma^2(\Omega)$.
%, the Stokes operator on $L_\sigma^p(\Omega)$ is defined as $A_p=-\mathcal{P}\Delta$ with domain
%Since $A_{p_1}$ and $A_{p_2}$ coincide on the intersection of their domains for $p_1$, $p_2\in(1,\infty)$, we will drop the index in the following.

Within this framework, by straightforward adaptation of arguments  in %according to results of
\cite{Lankeitnon} with only some necessary modifications, one can see that %it is known
the problem  \eqref{1.1}, \eqref{1.4}, \eqref{1.5} admits a global classical solution $(n,c,u,P)$  %in  two-dimensional settings
whenever $ \chi\in (0,1)$, $r\in \mathbb{R}$ and $\mu>0$,
which is unique up to addition of constants in the pressure variable $P$, and satisfies $n > 0 $,  $ c > 0$ in $\Omega\times [0,\infty)$.
 The first of our main results is concerned with the global boundedness of the solution %to \eqref{1.1}--\eqref{1.3}
 as well as its asymptotic behavior.

\begin{theorem}\label{theorem1.1} Let $f(n)=rn-\mu n^2$, $r\in \mathbb{R}$, $\mu>0$ and $\phi\in W^{1,\infty}(\Omega)$, and suppose that  $(n_0,c_0,u_0)$ satisfy \eqref{1.6}. %$\Omega\subset \mathbb{R}^2$  be a bounded domain with smooth boundary, and
If $(n,c,u,P)$ denotes the corresponding  global classical solution  to \eqref{1.1}, \eqref{1.4}, \eqref{1.5},
 % Assume that $ \chi\in (0,1)$ and $r\geq 0, \mu>0$,
then there exists  a  %an explicit
value
$\mu_0:=\mu_0(\Omega,\chi, r)\geq 0$  with $\mu_0(\Omega,\chi, 0)= 0 $ such that whenever $ \mu>\mu_0$, $(n,c,u)$ is global bounded, %and when $r>0$,
$$\|n(\cdot,t)-\frac {r_+} \mu\|_{L^\infty(\Omega)}\rightarrow 0, ~~
\|\frac{\nabla c}c(\cdot,t)\|_{L^\infty(\Omega)} \rightarrow 0,~~
\|u(\cdot,t)\|_{L^\infty(\Omega)} \rightarrow 0 ~~
$$
and when $r>0$,
$$~~\|c(\cdot,t) \|_{L^\infty(\Omega)} \rightarrow 0~~~~$$
as $ t \rightarrow \infty$.
%uniformly with respect to $x\in \Omega$.
%for,. there exists
\end{theorem}

As indicated in the above discussion, we need to introduce new ideas to show how
the regularizing effect of  the quadratic degradation in the chemotaxis--fluid model \eqref{1.1}
can counterbalance  the strongly destablizating action of chemotactic cross--diffusion caused by the combination of the consumption of $c$ with
the singular chemotaxis sensitivity of $n$. Specifically, we develop the conditional energy functional method in \cite{Winklerp} to show the global boundedness of solutions in the case of $r>0$, in which the key point is %, the most critical step under our approach
 to verify that % examine the behavior of the functionals
\begin{equation}\label{3.1}
\mathcal{F}(n,w):= \int_\Omega H(n) +\frac \chi2\displaystyle \int_\Omega |\nabla w|^2, ~~w:= -\ln (\frac c{\|c_0\|_{L^\infty(\Omega)}})
 \end{equation}
 with $H(s):=s\ln \frac  {\mu s}{er}+ \frac r \mu$ constitutes an energy functional in the sense that $\mathcal{F}(n,w)$ is non-increasing
  in time whenever $\mu$ is appropriately large relative to $r$ (see Lemma \ref{lemma32}). Indeed, from \eqref{3.13}, one can obtain the global bound of
$\int_\Omega n|\ln n| dx$ and $\int_\Omega |\nabla w|^2 dx$, which then serves as a starting point % the foundation
to  derive the uniform bound  of $\|n(\cdot,t)\|_{L^\infty(\Omega)}$
via the Neumann heat semigroup estimates. Furthermore, by making appropriate use of the dissipative information expressed in \eqref{3.13}, we can establish the convergence result asserted in Theorem  1.1. It is noted that compared to that of the case $r>0,\mu>0$, the proof of Theorem 1.1 in the case of $r\leq 0, \mu>0$ involves a more delicate analysis.
 In fact, unlike in the case $r=\mu=0$  or  $r>0,\mu>0$, \eqref{1.1} with $r\leq 0,\mu>0$ seems to lack  the favorable structure %in the case of $r\leq 0,\mu>0$,
that facilitates % accordingly
such conditional energy-type inequalities. %this  makes the %availability  of the decreasing conditional energy functional.
Taking full advantage of the decay information on $n$  in  $L^1-$norm expressed in \eqref{2.2}, our approach toward Theorem 1.1 is to  construct the quantity
 \begin{equation}\label{3.1}
\mathcal{F}(n,w):= \int_\Omega n(\ln n+a) +\frac \chi2\displaystyle \int_\Omega |\nabla w|^2
 \end{equation}
with parameter $a>0$ determined below (see \eqref{3.48}). Unlike in the case of $r>0$,  $\mathcal{F}(n,w)$
does not  enjoy monotonicity property, %$\mathcal{F}(n,w)$
it however satisfies a favorable  non-homogeneous differential inequality \eqref{3.55} in the sense that it can provide  us  %to derive
a  priori  information on solution such as  the global bound of
$\int_\Omega n|\ln n| dx $ and $\int_\Omega |\nabla w|^2 dx$ (see Lemma \ref{lemma34}), as well as
$\displaystyle\lim_{t\rightarrow \infty}\int_\Omega |\nabla w(\cdot,t)|^2=0$ (see \eqref{3.61}).

%as |?w|2 is suitably small  for  along trajectories of the boundary value problem in \eqref{2.1}.

As an important step to understand the  model \eqref{1.1} more comprehensively, we shall consider the %the asymptotic behavior and
 convergence rate of its classical solutions in the form of the following result:
%Beyond the boundedness of,

\begin{theorem}\label{theorem1.2} Let the assumptions of Theorem 1.1 hold and $r>0$.  Then one can find  $\mu_*(\chi, \Omega, r)>0$ such that %there exists $\theta^*\geq \theta_*$ such that
if $ \mu >\mu_*(\chi, \Omega, r)$,  the  classical solution of \eqref{1.1}, \eqref{1.4}, \eqref{1.5} presented in Theorem 1.1  satisfies
$$\|n(\cdot,t)-\frac {r} \mu\|_{L^\infty(\Omega)}\rightarrow 0, ~~
\|c(\cdot,t) \|_{L^\infty(\Omega)} \rightarrow 0,~~~~
\|u(\cdot,t) \|_{L^\infty(\Omega)} \rightarrow 0
$$
%in $(L^\infty(\Omega))^3$
as well as
$$\|\frac{\nabla c}c(\cdot,t)\|_{L^p(\Omega)} \rightarrow 0
$$
%in $L^\infty(\Omega)\times L^\infty(\Omega)$ and   $L^p(\Omega)$
for all $p>1$ exponentially as $t\rightarrow \infty$.
\end{theorem}
This implies that suitably large $\mu$ relative to $r$ enforces asymptotic stability of the corresponding constant equilibria of
\eqref{1.1}; however, the optimal lower bound on $\frac \mu r$ seems yet lacking.
The main ingredient of the  our approach toward Theorem 1.2  involves a so-called self-map-type reasoning. More precisely,
making use of the convergence properties of  $(n,\frac{|\nabla c|}c)$  asserted in Theorem 1.1, we prove by a self-map-type reasoning that whenever $\mu$ is suitably large compared with $r$,
$$(n(\cdot,t)-\frac r \mu, c(\cdot,t),u)\longrightarrow (0,0,0)
~~ \hbox{and}~~\frac{|\nabla c|}c(\cdot,t)\longrightarrow 0
$$
in $(L^\infty(\Omega))^3 $ and   $L^6(\Omega)$
exponentially as $t\rightarrow \infty$, respectively (see Lemma \ref{lemma43}).

As aforementioned, the limit case $r=0$ becomes relevant in several applications. In this limiting situation, the total cell population can readily be seen to decay in the large time limit
(cf. Lemma 2.3 below). As a consequence, we can obtain the decay properties of solutions, namely that the decay on $n$ in $L^1$  actually occurs   in $L^\infty$, and also for $c$. More precisely, our result reads as follows:

\begin{theorem}\label{theorem1.3} Let the assumptions of Theorem 1.1 hold and $r=0$. Then the classical solution  of \eqref{1.1}, \eqref{1.4}, \eqref{1.5} from  Theorem 1.1 satisfies
$$(n,c,\frac{|\nabla c|}c,u)\longrightarrow (0,0,0,0)
$$
in $(L^\infty(\Omega))^4$   algebraically  as $t\rightarrow \infty$.
\end{theorem}

The result indicates that structure generating dynamics in the spatially
two-dimensional version of \eqref{1.1},  \eqref{1.4},  \eqref{1.5}, if at all, occur on intermediate time
scales rather than
in the sense of a stable large-time pattern formation process.
Apparently, it leaves open  the
questions whether the more
colorful large time behavior can appear in the three-dimensional version  of  \eqref{1.1}. %ity of the physical setting
%can enforce more
%colorful large time behavior in \eqref{1.1}.

Our approach toward Theorem 1.3 uses an alternative method, which, at its core, is based on
the argument that the $L^\infty$-norm of $n$ can be controlled from above by appropriate multiples of $\frac 1{t+1}$.  This results from a suitable
 variation-of-constants representation of $n$, by which and in view of the decay information on $|\nabla w|$ in $L^\infty(\Omega)$,   %again by utilizing \eqref{2.2},
 the $L^1$ decay information on $u$ from \eqref{2.2} can be turned into
the $L^\infty$-norm of $n$ (see Lemma \ref{lemma45}). As a consequence, by comparison argument, we have a pointwise upper estimate for $w$ as well as a lower estimate for  $v$ (see Lemma \ref{lemma46}). Using $L^p-L^q$ estimates for the Neumann heat semigroup $(e^{t\Delta})_{t>0}$, we then successively show that  $\|\nabla w\|_{L^\infty}$
 and %$L^\infty(\Omega)$-norm of
$\|n\|_{L^\infty(\Omega)}$ can be controlled by  appropriate multiples of $\frac 1{t+1}$ from above and below, respectively (see Lemma \ref{lemma47}).
These  a priori estimates allow us to get the  pointwise lower estimate for $w$ as well as the upper estimate for  $c$, which complement the lower
bound for $c$ previously obtained,  and thereby prove that $c$ actually decays  algebraically.

\section{Preliminaries}  %Local existence and basic properties}
%\vspace{-1em}
In this section, we begin by recalling the  important $L^p-L^q$ estimates  for the Neumann heat semigroup $(e^{t\Delta})_{t>0}$  on
  bounded domains (see \cite{Cao, Winkler7}).
   %well-known Neumann heat semigroup estimates given in.  %mainly needed in the  proof of Theorem 1.2.
%We $L^p-L^q$ estimates for the Neumann heat semigroup.

\begin{lemma}\label{Lemma 2.1} %(Lemma 1.3 of \cite{Winkler7})
Let $(e^{t\Delta})_{t>0}$ denote the Neumann heat semigroup in the bounded domain $\Omega\subset\mathbb R^n$ and $\lambda_1>0$ denote the first nonzero eigenvalue of $-\Delta$ in $\Omega$ under the Neumann boundary condition. Then there exists $c_i>0$ ($i=1,2,3$) such that for all $t>0$,

%$(i)$ If $1\leq q\leq p\leq\infty$, then, for all $\omega\in L^q(\Omega)$ with $\int_\Omega\omega=0$,
%$$\|e^{t\Delta}\omega\|_{L^p(\Omega)}\leq c_1\left(1+t^{-\frac n2(\frac1q-\frac1p)}\right)e^{-\lambda_1t}\|\omega\|_{L^q(\Omega)};$$

$(i)$ If $1\leq q\leq p\leq\infty$, then for all $\omega\in L^q(\Omega)$,
$$\|\nabla e^{t\Delta}\omega\|_{L^p(\Omega)}\leq c_1\left(1+t^{-\frac12-\frac n2(\frac1q-\frac1p)}\right)e^{-\lambda_1t}\|\omega\|_{L^q(\Omega)};$$

$(ii)$ If $2\leq q\leq p<\infty$, then for all $\omega\in W^{1,q}(\Omega)$,
$$\|\nabla e^{t\Delta}\omega\|_{L^p(\Omega)}\leq c_2\left(1+t^{-\frac n2(\frac1q-\frac1p)}\right)e^{-\lambda_1t}\|\nabla\omega\|_{L^q(\Omega)};$$

$(iii)$ If $1\leq q\leq p\leq\infty$ or $1<q<\infty$ and $p=\infty$, then for all $\omega\in (L^q(\Omega))^n$,
$$\|e^{t\Delta}\nabla\cdot\omega\|_{L^p(\Omega)}\leq c_3\left(1+t^{-\frac12-\frac n2(\frac1q-\frac1p)}\right)e^{-\lambda_1t}\|\omega\|_{L^q(\Omega)}.$$
\end{lemma}\vspace{-1em}

%Let us also recall those special cases of the well-known Gagliardo--Nirenberg inequality
\begin{lemma}\label{Lemma2.2}
(Gagliardo--Nirenberg inequality) Let $ \Omega\subset \mathbb{R}^2$ be a bounded Lipschitz domain.
Then i) there is $ K_1> 0$ such that
$$ \|\nabla \varphi\|^4_{L^4(\Omega)}\leq K_1\|\Delta\varphi\|^2_{L^2(\Omega)}\|\nabla\varphi\|^2_{L^2(\Omega)}
$$for  all $\varphi\in W^{2,2}(\Omega)$  fulfilling  $\frac{\partial\varphi}{\partial \nu}|_{\partial\Omega}=0 $;

ii) there is  $ K_2> 0$ such that
$$ \|\varphi\|^3_{L^3(\Omega)}\leq K_2\| \varphi\|^2_{W^{1,2}(\Omega)}\|\varphi\|_{L^1(\Omega)}$$
for  all $\varphi\in W^{1,2}(\Omega)$.
\end{lemma}\vspace{-1em}

  In order to derive some essential estimates, it would be  more convenient to deal with a  nonsingular chemotaxis term of the form  $\nabla\cdotp (n \nabla w)$  instead  of $\nabla\cdotp (\frac {n}c\nabla c)$ in \eqref{1.1}. To this end, we employ the following transformation as  in \cite{Lankeitnon,Lankeitjde,Winkler3M}: %of $v$ %adapt take changing variables in \eqref{1.1} (see \cite{Lankeitjde,Lankeitnon,Winkler3M}), we
$$w:= -\ln (\frac c{\|c_0\|_{L^\infty(\Omega)}}), $$
whereupon $0\leq w\in C^0(\bar{\Omega}\times (0,\infty))\cap  C^{2,1}(\bar{\Omega}\times (0,\infty))$, and the problem \eqref{1.1}, \eqref{1.4}, \eqref{1.5} transforms to
\begin{equation}\label{2.1}
 \left\{\begin{array}{ll}
  n_t+u\cdot \nabla n =\triangle n+\chi \nabla\cdotp ( n\nabla w)+n(r- \mu n),& x\in \Omega, t>0,\\
 w_t+u\cdot \nabla w =\triangle w- |\nabla w|^2 +n, &x\in \Omega, t>0,\\
u_t+ (u\cdot \nabla) u=\Delta u+\nabla P+n\nabla\phi, & x\in \Omega, t>0,\\
\nabla \cdot u=0,& x\in \Omega, t>0,\\
 \nabla n\cdot \nu=\nabla w\cdot \nu=0, \quad u=0, &x\in \partial\Omega, t>0,\\

n(x,0)=n_0(x),~ w(x,0)=-\ln (\frac {c_0(x)}{\|c_0\|_{L^\infty(\Omega)}}), ~~ u=u_0(x), & x\in \Omega.
   \end{array}\right.
\end{equation}

Let us first recall
some basic but important information about  $(n,w)$  due to the presence of the quadratic  degradation term in the first equation of \eqref{2.1}.
\begin{lemma}\label{lemma23}
The classical solution $(n,w,u,P)$ of \dref{2.1} satisfies

(i) $\displaystyle\lim\sup_{t\rightarrow \infty}\|n(\cdot,t)\|_{L^1(\Omega)}\leq\frac {|\Omega|r_+} \mu;$

(ii)$\displaystyle \int^{t}_{t_0} \|n(\cdot,s)\|^2_{L^2(\Omega)}ds\leq  \frac {r_+} \mu \int^{t}_{t_0}\|n(\cdot,s)\|_{L^1(\Omega)}ds
+   \frac 1 \mu  \|n(\cdot,t_0)\|_{L^1(\Omega)}$  for all  $t> t_0$;

(iii) $\displaystyle  \int^{t}_{t_0}  \int_\Omega |\nabla w|^2 dxds\leq \int_\Omega w(x,t_0) dx
+\int^{t}_{t_0}\|n(\cdot,s)\|_{L^1(\Omega)}ds$  for all  $t>t_0$.

In particular, if $r\leq 0$, then
\begin{equation}\label{2.2}
\|n(\cdot,t)\|_{L^1(\Omega)}\leq\displaystyle\frac {|\Omega|} {\mu(t+\gamma)}~~~\mbox{for all}~~ t>t_0\end{equation}
 with $\gamma=\displaystyle\frac {|\Omega|}{\mu\int_\Omega n_0(x) dx}$.

\end{lemma}
\it Proof.  \rm Integrating the first equation in \eqref{2.1} and using %According to \eqref{2.1} and
the Cauchy--Schwarz inequality, we get
\begin{equation}\label{2.3}
\begin{array}{rl}
\displaystyle\frac d{dt}\int_\Omega n= &\displaystyle r\int_\Omega n- \mu \int_\Omega n^2 \leq  \displaystyle r_+\int_\Omega n - \frac\mu {|\Omega|}(\int_\Omega n )^2
 \end{array}
 \end{equation}
which yields (i) readily.  By the time integration of \eqref{2.3} over $(t_0,t)$, we get (ii) immediately.
In addition, from the second equation in \eqref{2.1}, $\nabla\cdot u=0$ and $u=0$ on  $\partial\Omega$, it follows that
\begin{equation}\label{2.4}
\displaystyle\frac d{dt}\int_\Omega w = -\displaystyle \int_\Omega |\nabla w|^2 +\int_\Omega n,
 \end{equation}
and thus establishes  (iii).

When $r\leq 0$, it follows from \eqref{2.3} that
\begin{equation}\label{2.4}
\displaystyle\frac d{dt}\int_\Omega n \leq  \displaystyle - \frac\mu {|\Omega|}(\int_\Omega n )^2
 \end{equation}
which then yields \eqref{2.2} by the time integration.

In order to make use of the spatio-temporal properties provided by Lemma 2.3(ii) to estimate the ultimate  bound  of $\int_\Omega |\nabla u|^2$, we shall utilize the following elementary lemma (see Lemma 3.4 of \cite{WinklerJFA}):

\begin{lemma}\label{lemma24d}Let $t_0\geq 0, T\in (t_0,\infty]$, $a>0$ and $b>0$,  and suppose that the nonnegative function  $h\in L^1_{loc}(\mathbb{R})$ satisfies $\int^{t+1}_{t} h(s)ds\leq b$ for all $t\in [t_0,T]$. If $y \in C^0([t_0, T ))\cap C^1([t_0, T )) $
has the property that
$$
y'(t) + ay(t) \leq  h(t)~~\hbox{for all}~~ t\in (t_0, T ),
$$
then
$$
y(t)\leq e^{-a(t-t_0)}y(t_0)+\frac b{1-e^{-a}} ~~\hbox{for all }~~t\in [t_0, T ).
$$

\end{lemma}

With
%In light of
Lemmas \ref{lemma23} and \ref{lemma24d} at hand, we can employ the standard energy inequality associated with the fluid
evolution system  in  \dref{2.1} to  derive some  boundedness results for $u$.
%of the natural energy functional %in a standard manner
 %estimate  $ \| u(\cdot,t)\|_{W^{1,2}(\Omega)}$.  We first obtain the following bounds of $\| u(\cdot,t)\|_{L^2(\Omega)}$.
\begin{lemma}\label{lemma25}
For the global  classical solution $(n,w,u)$ of \dref{2.1}, we have
% There exists some $t_0>0$ such that for all $t>t_0$,

i) if $r>0$, then
\begin{align}\label{2.6}
 \displaystyle\lim\sup_{t\rightarrow \infty}\|u(\cdot,t)\|^2_{L^2(\Omega)}\leq \frac {3(1+r) |\Omega|} \mu
\frac {\|\nabla\phi \|^2_{L^\infty (\Omega)}}{C_p(1-e^{-\frac {C_p}2})}\frac {r} \mu
 % \|u(\cdot,t_0)\|_{L^2(\Omega)}e^{-\frac {C_p}2(t-t_0)}+\frac b{1-e^{-\frac {C_p}2}}.
\end{align}
as well as
\begin{align}\label{2.7}
\displaystyle\lim\sup_{t\rightarrow \infty}\int^{t+1}_{t}\|\nabla u(\cdot,s)\|^2_{L^2(\Omega)}ds\leq \frac {5(1+r) |\Omega|} \mu
\frac {\|\nabla\phi \|^2_{L^\infty (\Omega)}}{C_p(1-e^{-\frac {C_p}2})}\frac {r} \mu
\end{align}
with Poincar\'{e} constant $C_P>0$.

ii) if $r\leq 0$, then
\begin{align}\label{2.8}
\int_\Omega|u(\cdot,t)|^2\leq  \|u(\cdot,t_0)\|^2_{L^2(\Omega)}e^{-\frac {C_p}2(t-t_0)}+
\frac {2 |\Omega|} {\mu^2}
\frac {\|\nabla\phi \|^2_{L^\infty (\Omega)}}{C_p(1-e^{-\frac {C_p}2})}\frac {1} {t_0+\gamma}
~~~~~\hbox{for all}~~t>t_0
\end{align}
as well as
\begin{align}\label{2.9}
& \int^{t+1}_{t}\|\nabla u(\cdot,s)\|^2_{L^2(\Omega)}ds\nonumber\\
\leq &
\|u(\cdot,t_0)\|^2_{L^2(\Omega)}e^{-\frac {C_p}2(t-t_0)}+
\frac {4 |\Omega|} {\mu^2}
\frac {\|\nabla\phi \|^2_{L^\infty (\Omega)}}{C_p(1-e^{-\frac {C_p}2})}\frac {1} {t_0+\gamma}
~~~\hbox{for all}~t>t_0.
\end{align}
\end{lemma}
\it Proof. \rm
i)  According to the Poincar\'{e} inequality, one can find some constant $C_p>0$ such that
$$
C_p\int_\Omega|u|^2\leq \int_\Omega|\nabla u|^2.
 $$

 Testing the third equation in \dref{2.1}  by $u$  and using the   H\"{o}lder inequality,  we obtain
\begin{align*}
\frac{d}{dt}\int_\Omega|u|^2+C_p\int_\Omega|u|^2+\int_\Omega|\nabla u|^2&\leq 2\int_\Omega n\nabla\phi \cdot u
\\
&\leq 2\|\nabla\phi \|_{L^\infty (\Omega)}\|n\|_{L^2(\Omega)}\|u\|_{L^2(\Omega)}
\\
&\leq \frac {C_p}2\| u\|_{L^2(\Omega)}^2+\frac 2{C_p}\|\nabla\phi \|^2_{L^\infty (\Omega)}\|n\|_{L^2(\Omega)}^2,
\end{align*}
due to $u|_{\partial \Omega} = 0$ and  $\nabla \cdot u = 0$.

Writing $h(t)=\frac 2{C_p}\|\nabla\phi \|^2_{L^\infty (\Omega)}\|n(\cdot,t)\|_{L^2(\Omega)}^2$, we see that $y(t):=\int_\Omega|u(\cdot,t)|^2$
satisfies
\begin{align}\label{2.10}
y'(t) +\frac {C_p}2 y(t)+ \int_\Omega|\nabla u(\cdot,t)|^2 \leq h(t)~~~~~\hbox{for all}~~t>0.
\end{align}

In view of Lemma \ref{lemma23} (i) and (ii), we know that
\begin{align}\label{2.11}
\displaystyle\lim\sup_{t\rightarrow \infty}\int^{t+1}_{t} h(s)ds\leq \frac 2{C_p}\|\nabla\phi \|^2_{L^\infty (\Omega)}
\frac {(1+r) |\Omega|} \mu\frac {r} \mu.
\end{align}
An application of Lemma \ref{lemma24d} thus shows that there exists positive $t_0>0$ such that
$$\int_\Omega|u(\cdot,t)|^2\leq  \|u(\cdot,t_0)\|^2_{L^2(\Omega)}e^{-\frac {C_p}2(t-t_0)}+
\frac {3(1+r) |\Omega|} \mu
\frac {\|\nabla\phi \|^2_{L^\infty (\Omega)}}{C_p(1-e^{-\frac {C_p}2})}\frac {r} \mu
~~~~~\hbox{for all}~~t>t_0
$$
and thereby verifies \eqref{2.6}. Thereafter, again thanks to \eqref{2.11}, an integration of  \eqref{2.10} in time yields
\eqref{2.7}.

ii) In view of \eqref{2.2}, we have
\begin{align}\label{2.12}
\int^{t+1}_{t} h(s)ds\leq \frac 2{C_p}\|\nabla\phi \|^2_{L^\infty (\Omega)}
\frac { |\Omega|} {\mu^2}
\frac {1} {t+\gamma},
\end{align}
whereupon Lemma \ref{lemma24d} guarantees that
$$\int_\Omega|u(\cdot,t)|^2\leq  \|u(\cdot,t_0)\|^2_{L^2(\Omega)}e^{-\frac {C_p}2(t-t_0)}+
\frac {2 |\Omega|} {\mu^2}
\frac {\|\nabla\phi \|^2_{L^\infty (\Omega)}}{C_p(1-e^{-\frac {C_p}2})}\frac {1} {t_0+\gamma}
~~~~~\hbox{for all}~~t>t_0.
$$
This  precisely warrants \eqref{2.8}, and thereby in turn  yields \eqref{2.9}
after integrating  \eqref{2.10} over $(t,t+1)$ and once more employing \eqref{2.12}.

 Now by a further testing procedure, we can turn the above information into the estimate of $\|\nabla u(\cdot,t)\|_{L^2(\Omega)}$,
 particularly its decay in the case of $r=0$, %the spatial gradient of $u $
on the basis of an interpolation argument, which is inspired by an approach illustrated in section 3.2 of \cite{TaoZAMP}.
% which essentially relies on our assumption that the spatial setting is two-dimensional.

 \begin{lemma}\label{lemma26} For the global  classical solution $(n,w,u,P)$ of \dref{2.1}, we have

i) if $r>0$, then  there exists $\mu_1:=\mu_1(\Omega, r)>0$ such that for all $\mu>\mu_1$,
\begin{align}\label{2.13}
\displaystyle\lim\sup_{t\rightarrow \infty}\|\nabla u(\cdot,t)\|_{L^2(\Omega)}\leq \frac {1}{17K_1|\Omega|}
%\frac {2(1+r) |\Omega|} \mu
%\frac {\|\nabla\phi \|_{L^\infty (\Omega)}(2-e^{-\frac {C_p}2})}{C_p(1-e^{-\frac {C_p}2})}\frac {r} \mu.
 \end{align}

ii) if $r\leq 0$, then for any $\mu>0$,
\begin{align}\label{2.14}
\displaystyle\lim_{t\rightarrow \infty}\|\nabla u(\cdot,t)\|_{L^2(\Omega)}=0.
 \end{align}

\end{lemma}
\it Proof. \rm  Applying the Helmholtz projector $\mathcal{P}$ to the third equation in \eqref{2.1},  multiplying the
resulting identity
$ u_t + Au = -\mathcal{P}[(u \cdot\nabla)u] +\mathcal{P}[n\nabla \phi]$
 by $A u$, and using the Gagliardo--Nirenberg inequality,  we can find $C_1>0$ such that
    \begin{align*}
\frac12\frac{d}{dt}\int_{\Omega}|\nabla u|^2+\int_\Omega|A u|^2&=-\int_\Omega \mathcal{P} [(u \cdot \nabla)u] \cdot A u+\int_\Omega \mathcal{P}[n \nabla\phi] \cdot A u
\\
&\leq \frac 12\int_\Omega|A u|^2+\int_\Omega |(u \cdot \nabla)u|^2+\|\nabla\phi \|^2_{L^\infty (\Omega)}\int_\Omega n^2\\
&\leq \frac12\int_\Omega|A u|^2+\| u\|^2_{L^\infty(\Omega)} \|\nabla u\|^2_{L^2(\Omega)}+\|\nabla\phi \|^2_{L^\infty (\Omega)}\int_\Omega n^2\\
&\leq \frac12\int_\Omega|A u|^2+C_1 \|A u\|_{L^2(\Omega)}\|u\|_{L^2(\Omega)} \|\nabla u\|^2_{L^2(\Omega)}+
\|\nabla\phi \|^2_{L^\infty (\Omega)}\int_\Omega n^2\\
&\leq \int_\Omega|A u|^2+\frac{C^2_1}2 \|u\|^2_{L^2(\Omega)} \|\nabla u\|^4_{L^2(\Omega)}+
\|\nabla\phi \|^2_{L^\infty (\Omega)}\int_\Omega n^2,
\end{align*}
which entails   $y(t):=\int_\Omega|\nabla u(\cdot,t)|^2$
satisfies
\begin{align}\label{2.15}
y'(t) \leq h_1(t)y(t)+h_2(t)~~~~~\hbox{for all}~~t>0
\end{align}
with $h_1(t)=C^2_1\|u(\cdot,t)\|^2_{L^2(\Omega)} \|\nabla u(\cdot,t)\|^2_{L^2(\Omega)}$ and $h_2(t)=2\|\nabla\phi \|^2_{L^\infty (\Omega)}\| n (\cdot,t)\|^2_{L^2(\Omega)}$.

i) In order to prepare the  integration of \eqref{2.15}, we may use Lemma \ref{lemma25} i) to find some $t_0>0$ such that
\begin{align*}
\|u(\cdot,t)\|^2_{L^2(\Omega)}\leq C_2:=\frac {3(1+r) |\Omega|} \mu
\frac {\|\nabla\phi \|^2_{L^\infty (\Omega)}}{C_p(1-e^{-\frac {C_p}2})}\frac {r} \mu
\end{align*}
and
\begin{align*}
\int^{t}_{t-1}\|\nabla u(\cdot,s)\|^2_{L^2(\Omega)}ds\leq 2 C_2
\end{align*}
 for all $t>t_0+1$.

Hence for any $t>t_0+1$, we can find $t_*=t_*(t)\in [t-1,t)$ such that
\begin{align}\label{2.16}
\|\nabla u(\cdot,t_*)\|^2_{L^2(\Omega)}\leq 2 C_2,
\end{align}
and then integrating \eqref{2.15} over $(t_*,t)$ yields
\begin{align*}
y(t)&\leq y(t_*)e^{\int ^t_{t_*}h_1(\sigma)d\sigma} +  \int^t_{t_*} e^{\int ^t_{s}h_1(\sigma)d\sigma} h_2(s)ds\\
&\leq (2+C_p)C_2e^{2C_1^2C_2^2}
\end{align*}
and thereby verifies \eqref{2.13}.

ii) For any $t_0>1$ and $t>t_0+2$, we use Lemma \ref{lemma25} ii) to pick $t_*=t_*(t)\in [t-1,t)$ fulfilling
\begin{align*}\label{2.19}
\|\nabla u(\cdot,t_*)\|^2_{L^2(\Omega)}=&
\int^{t}_{t-1}\|\nabla u(\cdot,s)\|^2_{L^2(\Omega)}ds\\
\leq & \|u(\cdot,t_0)\|^2_{L^2(\Omega)}e^{-\frac {C_p}2(t-1-t_0)}+
\frac {4 |\Omega|} {\mu^2}
\frac {\|\nabla\phi \|^2_{L^\infty (\Omega)}}{C_p(1-e^{-\frac {C_p}2})}\frac {1} {t_0+\gamma},
\end{align*}
as well as
\begin{align*}
\int ^t_{t-1}h_1(\sigma)d\sigma &\leq C^2_1\max_{t-1\leq s\leq t}\|u(\cdot,s)\|^2_{L^2(\Omega)}\int ^t_{t-1} \|\nabla u(\cdot,s)\|^2_{L^2(\Omega)}ds\\
 &\leq C^2_1 (\|u(\cdot,t_0)\|^2_{L^2(\Omega)}e^{-\frac {C_p}2(t-1-t_0)}+
\frac {4 |\Omega|} {\mu^2}
\frac {\|\nabla\phi \|^2_{L^\infty (\Omega)}}{C_p(1-e^{-\frac {C_p}2})}\frac {1} {t_0+\gamma})^2.
\end{align*}
In  addition, by \eqref{2.12} we also have
 \begin{align*}
\int ^t_{t-1}h_2(\sigma)d\sigma = 2\|\nabla\phi \|^2_{L^\infty (\Omega)}\int ^t_{t-1}\| n (\cdot,s)\|^2_{L^2(\Omega)}ds
\leq
 2\|\nabla\phi \|^2_{L^\infty (\Omega)}
\frac { |\Omega|} {\mu^2}
\frac {1} {t-1+\gamma}.\end{align*}

Therefore combining the above inequalities, \eqref{2.15} implies that
\begin{align*}
y(t)&\leq y(t_*)e^{\int ^t_{t-1}h_1(\sigma)d\sigma} + e^{\int ^t_{t-1}h_1(\sigma)d\sigma} \int^t_{t-1} h_2(s)ds
\end{align*}
and thus  \eqref{2.14} holds readily.

\section{Global boundedness of solutions}
In this section, we show that  the classical solution of the problem  \eqref{2.1} is globally  bounded
in the cases of $r>0$ and $r\leq  0$, respectively.

\subsection{The case $r>0$}
In this subsection, we derive the global boundedness of  solutions to \eqref{2.1}  whenever $\mu$ is suitably large compared with $r$. As in \cite{Winklerp}, the  main idea is to examine the behavior of the functional
\begin{equation}\label{3.1}
\mathcal{F}(n,w):= \int_\Omega H(n) +\frac \chi2\displaystyle \int_\Omega |\nabla w|^2
 \end{equation}
 where $H(s):=s\ln \frac  {\mu s}{er}+ \frac r \mu$, along trajectories of the boundary value problem \eqref{2.1}.

The following elementary property of $H(n)$ will be used in the sequel.
\begin{lemma}\label{lemma31} For all nonnegative function $n\in C(\bar{\Omega})$,
$H(n)\geq 0$.
\end{lemma}
\it Proof.  \rm It is easy to verify that $H(\frac r \mu)=0$,  $H'(\frac r \mu)=0$ and $H''(s)=\frac 1{s}\geq 0$, which implies $H(n)\geq 0$ for all $n\geq 0$.

Now we can describe the evolution  of $\mathcal{F}(n,w)$  along the trajectories of (2.1)  by the standard testing procedure.
\begin{lemma}\label{lemma32}
Let $\Omega\subset \mathbb{R}^2$ be a smooth bounded domain and $(n,w,u)$ be the global classical solution of \eqref{2.1} with $r>0, \mu>0$. Then there exists %$\mu_*>0$ and
$t_*>0$ such that
\begin{equation}\label{3.2}
\frac d{dt}{\mathcal{F}(n,w)}\leq 0 ~~\hbox{for all}~~ t\geq t_*
\end{equation}
whenever $\mu>\mu_2(\Omega,\chi,r):=\max\{\mu_1,\frac{K_1(36+16\chi)|\Omega|}{ \chi}r\} $.
\end{lemma}
\it Proof.  \rm  Multiplying  the first equation in (2.1) by $H(n)$ and integrating by parts,  we get
\begin{equation}\label{3.3}
\begin{array}{rl}
\displaystyle\frac{d}{dt}\int_\Omega H(n)=&\displaystyle\int_\Omega H'(n)(\triangle n+\chi \nabla\cdotp ( n\nabla w)+rn- \mu n^2-u\cdot \nabla n)\\
=&-\displaystyle\int_\Omega H''(n) ( |\nabla n|^2 +\chi n \nabla n\cdot\nabla w)+ \displaystyle\int_\Omega  H'(n)(rn- \mu n^2)\\[2mm]
=& -\displaystyle \displaystyle\int_\Omega  \frac{ |\nabla n|^2}n-\chi\displaystyle\int_\Omega \nabla n\cdot\nabla w
+ \displaystyle\displaystyle\int_\Omega  (\ln n-\ln\frac r\mu) (rn- \mu n^2)\\[2mm]
\leq & -\displaystyle \displaystyle\int_\Omega  \frac{ |\nabla n|^2}n-\chi\displaystyle\int_\Omega \nabla n\cdot\nabla w \\
\end{array}
\end{equation}
due to  $(\ln n-\ln\frac r\mu) (rn- \mu n^2)\leq 0$, $\nabla\cdot u=0$ and $u=0 $ on $\partial\Omega$.

On the other hand, testing the second equation in \eqref{2.1} by $-\triangle w$, using $\nabla\cdot u=0$ and $u=0 $ on $\partial\Omega$ again,
%.the Young inequality,
we can obtain
\begin{equation*}\label{3.4}
\begin{array}{rl}
\displaystyle\frac 1 2\frac{d}{dt} \int_\Omega |\nabla w|^2+  \int_\Omega |\triangle w|^2=&
\displaystyle
\int_\Omega |\nabla w|^2 \triangle w+ \displaystyle\int_\Omega \nabla n\cdot\nabla w+\int_\Omega (u\cdot\nabla w) \triangle w \\[2mm]
\leq &  \displaystyle\frac 1 2 \int_\Omega |\triangle w|^2+ \displaystyle\frac {1} 2 \int_\Omega |\nabla w|^4
 +\displaystyle\int_\Omega \nabla n\cdot\nabla w+\int_\Omega (u\cdot\nabla w) \triangle w \\[2mm]
 = &  \displaystyle\frac 1 2 \int_\Omega |\triangle w|^2+ \displaystyle\frac {1} 2 \int_\Omega |\nabla w|^4
 +\displaystyle\int_\Omega \nabla n\cdot\nabla w-\int_\Omega \nabla w\cdot (\nabla u\cdot\nabla w) .
\end{array}
\end{equation*}
Furthermore, by Lemma \ref{Lemma2.2} i) and the Cauchy--Schwarz inequality, we get
\begin{equation*}\label{3.4}
\begin{array}{rl}
\displaystyle\frac 1 2\frac{d}{dt} \int_\Omega |\nabla w|^2+ \displaystyle\frac 1 2 \int_\Omega |\triangle w|^2
\leq &  \displaystyle\frac {K_1} 2\|\nabla w\|^2_{L^2(\Omega)} \int_\Omega |\triangle w|^2
 +\displaystyle\int_\Omega \nabla n\cdot\nabla w+\displaystyle\int_\Omega |\nabla u||\nabla w|^2\\[2mm]
 \leq & ( \displaystyle\frac {K_1} 2\|\nabla w\|^2_{L^2(\Omega)}
 +K_1|\Omega|^{\frac12}\|\nabla u\|_{L^2(\Omega)})
  \int_\Omega |\triangle w|^2
 +\displaystyle\int_\Omega \nabla n\cdot\nabla w

\end{array}
\end{equation*}
and thus
\begin{equation}\label{3.4}
\displaystyle\frac 1 2\frac{d}{dt} \int_\Omega |\nabla w|^2+ \displaystyle\frac 1 2(1-K_1\|\nabla w\|^2_{L^2(\Omega)} -2
K_1|\Omega|^{\frac12}\|\nabla u\|_{L^2(\Omega)}) \int_\Omega |\triangle w|^2
\leq
 \displaystyle\int_\Omega \nabla n\cdot\nabla w.
\end{equation}

Since $2{\mathcal{F}(n,w)}\geq \chi\|\nabla w\|^2_{L^2(\Omega)}$ due to $H(n)\geq 0$, combining \eqref{3.4} with \eqref{3.3} yields
\begin{equation}\label{3.5}
\frac d{dt}{\mathcal{F}(n,w)} +\displaystyle \displaystyle\int_\Omega  \frac{ |\nabla n|^2}n+
(\displaystyle\frac \chi 2-  K_1{\mathcal{F}(n,w)}-2\chi
K_1|\Omega|^{\frac12}\|\nabla u\|_{L^2(\Omega)} ) \int_\Omega |\triangle w|^2\leq 0 ~~\hbox{for}~ t>0.
\end{equation}

On the other hand, from \eqref{2.13}, it is possible to pick some  $t_0>0$ such that
$$ 16
K_1|\Omega|^{\frac12}\|\nabla u(\cdot,t)\|_{L^2(\Omega)}<1 ~~~\hbox{for all }~ t>t_0
$$
whenever $\mu>\mu_1 $, and thereby
\begin{equation}\label{3.5b}
\frac d{dt}{\mathcal{F}(n,w)} +\displaystyle \displaystyle\int_\Omega  \frac{ |\nabla n|^2}n+
(\displaystyle\frac {3\chi} 8-  K_1\mathcal{F}(n,w)) \int_\Omega |\triangle w|^2
\leq 0 ~~\hbox{for }~ t>t_0.
\end{equation}

In what follows, we shall show that there exists $t_*>t_0 $  such that
$ 4 K_1\mathcal{F}(n,w)(t_*)<\chi$.

Firstly by Lemma \ref{lemma23}(i), there exists $t_1>t_0 $ such that for all $t>t_1$
\begin{equation}\label{3.6}
\|n(\cdot,t)\|_{L^1(\Omega)}\leq\frac {3|\Omega|r} {2\mu},
\end{equation}
which along with  Lemma \ref{lemma23}(iii) yields
\begin{equation*}
\begin{array}{rl}
\displaystyle  \int^{t_2}_{t_1}  \int_\Omega |\nabla w|^2 \leq
&\displaystyle\int_\Omega w(\cdot,t_1)
+\int^{t_2}_{t_1}\|n(\cdot,s)\|_{L^1(\Omega)}ds\\
\leq
&\displaystyle\int_\Omega w_0(x)
+\int^{t_1}_{0}\|n(\cdot,s)\|_{L^1(\Omega)}ds+\frac {3|\Omega|r} {2\mu}(t_2-t_1).
\end{array}
\end{equation*}
Similarly invoking Lemma 2.3(i) and (ii),  we find that
$$\displaystyle \int^{t_2}_{t_1} \|n(\cdot,s)\|^2_{L^2(\Omega)}ds\leq  \frac {3|\Omega|} {2}(\frac r \mu)^2 (t_2-t_1)
+   \frac 1 \mu  \|n(\cdot,t_1)\|_{L^1(\Omega)}.
$$
Hence there exists $t^*>t_1$ suitably large such that whenever $t_2\geq t^*$,
\begin{equation}\label{3.7}
\displaystyle  \int^{t_2}_{t_1}  \int_\Omega |\nabla w|^2 \leq
\frac {2|\Omega|r} \mu(t_2-t_1)
\end{equation}
and
\begin{equation}\label{3.8}
\displaystyle  \int^{t_2}_{t_1}  \|n(\cdot,s)\|^2_{L^2(\Omega)}ds \leq
2|\Omega| (\frac r \mu)^2 (t_2-t_1).
\end{equation}
Let
 $$\mathcal{S}_1:=\{t\in[t_1,t_2]| \int_\Omega |\nabla w(\cdot,t)|^2\geq\frac {8|\Omega|r} \mu \}
 $$
 and
$$\mathcal{S}_2:=\{t\in[t_1,t_2]| \|n(\cdot,t)\|^2_{L^2(\Omega)}\geq 8|\Omega| (\frac r \mu)^2 \}.
$$
 Then \begin{equation}\label{3.9}
 |\mathcal{S}_1|\leq \frac {|t_2-t_1|}4, ~~~|\mathcal{S}_2|\leq \frac {|t_2-t_1|}4.
 \end{equation}
  In order to estimate the size of
$ \mathcal{S}_1$ and $\mathcal{S}_2$,  we recall \eqref{3.7} to get
\begin{equation*}
\frac {8|\Omega|r} \mu|\mathcal{S}_1|\leq
\displaystyle  \int^{t_2}_{t_1}  \int_\Omega |\nabla w|^2 \leq
\frac {2|\Omega|r} \mu(t_2-t_1)
\end{equation*}
and thus  $|\mathcal{S}_1|\leq \frac {|t_2-t_1|}4$ is valid. Similarly, one can verify that $|\mathcal{S}_2|\leq \frac {|t_2-t_1|}4$.

As \eqref{3.9} warrants that
$$
|(t_1,t_2)\setminus (\mathcal{S}_1\cup  \mathcal{S}_2)| \geq  \frac {|t_2-t_1|}2,
$$  one can conclude that  there exists $t_*\in (t_1,t_2) $  such that
\begin{equation}\label{3.10}
\displaystyle \|n(\cdot,t_*)\|^2_{L^2(\Omega)}<  8|\Omega| (\frac r \mu)^2
\end{equation}
and
\begin{equation}\label{3.11}
\int_\Omega |\nabla w(\cdot,t_*)|^2<\frac {8|\Omega|r} \mu.
\end{equation}
Applying  $ \xi \ln \frac \xi \sigma\leq  \eta \xi^2+ \ln  \frac 1{\eta \sigma}\cdot \xi$ for all $\xi>0, \eta>0, \sigma >0$ %with $\sigma =\frac {er} \mu$
(see Lemma 5.5 of
\cite{Winklerp}),  and combining \eqref{3.6}  with \eqref{3.10}, we then arrive at
\begin{equation*}\begin{array}{rl}
\displaystyle\int_\Omega H(n)(\cdot,t_*) \leq& \displaystyle\frac {\mu} { r}\int_\Omega n^2(\cdot,t_*)
-\int_\Omega n(\cdot,t_*) +\frac {r}\mu|\Omega|
 \\ [3mm]
\leq &
\displaystyle \frac {9 |\Omega|r} \mu.  %+\frac {3}{\chi}|\Omega|\frac {r}{\mu}- \frac {3} {\chi}|\Omega| \frac {r}{\mu}\ln \frac r{\mu}.
\end{array}\end{equation*}
Thereupon  from \eqref{3.11} and   the definition of $ \mathcal{F}(n,w)$, it follows that
$$ \mathcal{F}(n,w)(t_*)< (9+ 4\chi)|\Omega|\frac r \mu,
$$
which entails that
$ 4 K_1\mathcal{F}(n,w)(t_*)<\chi$ provided $\mu> \frac{K_1(36+16\chi)|\Omega|r}{ \chi} $.

As an immediate consequence of \eqref{3.5}, we have
\begin{equation}\label{3.13}
\frac d{dt}{\mathcal{F}(n,w)} +\displaystyle \displaystyle\int_\Omega  \frac{ |\nabla n|^2}n+
\displaystyle\frac \chi 8 \int_\Omega |\triangle w|^2\leq 0 ~~\hbox{for all}~ t>t_*
\end{equation}
when $\mu>\mu_2(\Omega,\chi,r)$, and thus end the proof of this lemma. %is complete.

Additionally from \eqref{3.13}, one can  also conclude that

\begin{corollary}\label{Cor1}  %Let $(p,c,w)$ be the global classical solution of (1.1). Then
Under the conditions of Lemma 3.2, we have
\begin{equation}\label{3.14}
\mathcal{F}(n,w)(t) +\displaystyle\displaystyle\int^t_{t_*}\int_\Omega  \frac{ |\nabla n|^2}n+
\displaystyle\frac \chi 8 \int^t_{t_*}\int_\Omega |\triangle w|^2\leq (9+ 4 \chi)|\Omega|\frac r \mu~~~\hbox{for all}~ t>t_*.
\end{equation}

\end{corollary}

Next  by a further testing  procedure, we can turn the above information  into the uniform-in-time boundedness of $\|n(\cdot,t)\|_{L^2(\Omega)}$ and $\|\nabla w(\cdot,t)\|_{L^4(\Omega)}$ if $\mu$ is appropriately large compared with $r$, which will serve as the foundation for the proof of the global  boundedness of $\|n(\cdot,t)\|_{L^\infty(\Omega)}$ and  $\|\nabla w(\cdot,t)\|_{L^\infty(\Omega)}$.

\begin{lemma}\label{lemma33}
Under the assumptions in Lemma 3.2, there exists  $C>0$  % independent of  parameter $r>0$,
such that
\begin{equation}\label{3.15}
\displaystyle  \|n(\cdot,t)\|_{L^2(\Omega)}+\|\nabla w(\cdot,t)\|_{L^4(\Omega)}\leq   C~~\hbox{for all}~~ t\geq t_*
\end{equation}
provided  $\mu> \frac{182 K_2|\Omega| r } {\chi^2}$.
\end{lemma}
\it{Proof.}\rm\quad Multiplying the first equation in \eqref{2.1} by $n$ and integrating the result over $\Omega$, we get
\begin{equation}\label{3.16}
\begin{array}{rl}
\displaystyle  \frac12 \frac d{dt}\int_\Omega n^2= &
-\displaystyle\int_\Omega |\nabla n|^2- \chi\int_\Omega n \nabla n \nabla w +r \int_\Omega n^2- \mu \int_\Omega n^3\\[2mm]
\leq & -\displaystyle\frac12\displaystyle\int_\Omega |\nabla n|^2+  \frac12 \int_\Omega n^2 | \nabla w|^2 +r \int_\Omega n^2- \mu \int_\Omega n^3.
\end{array}
\end{equation}

On the other hand, by the second  equation  in \eqref{2.1} and the identity $\nabla w\cdot \nabla \Delta w=\frac 12 \Delta |\nabla w|^2- |D^2 w|^2$, we obtain
\begin{equation}\label{3.17}
\begin{array}{rl}
&\displaystyle \frac{d}{dt}\int_{\Omega}|\nabla w|^4\\
=
&\displaystyle  2\int_{\Omega}|\nabla w|^2 \triangle  |\nabla w|^2-
4 \int_{\Omega}|\nabla w|^2  |D^2 w|^2-
4 \int_{\Omega}|\nabla w|^2 \nabla w \cdot\nabla |\nabla w |^2\\[2mm]
&+
4 \displaystyle \int_{\Omega}|\nabla w|^2 \nabla n \cdot\nabla w
- 4 \displaystyle \int_{\Omega}|\nabla w|^2
\nabla w\cdot \nabla ( u \cdot\nabla w)\\[2mm]
= & - \displaystyle  2\int_{\Omega}|\nabla |\nabla w|^2|^2-
4 \int_{\Omega}|\nabla w|^2  |D^2 w|^2-
4 \int_{\Omega}|\nabla w|^2 \nabla w \cdot\nabla |\nabla w |^2-4\displaystyle  \int_{\Omega}n |\nabla w|^2   \triangle w
\\[3mm]
& -4\displaystyle  \int_{\Omega}n \nabla|\nabla w|^2\cdot  \nabla w +
2 \displaystyle  \int_{\partial\Omega} |\nabla w|^2 \frac{ \partial|\nabla w|^2}{\partial\nu}- 4 \displaystyle \int_{\Omega}|\nabla w|^2
\nabla w\cdot (\nabla u \cdot\nabla w)
\end{array}
\end{equation}
due to $\nabla\cdot u=0$ and $u=0 $ on $\partial\Omega$.

According to
$$\frac{ \partial|\nabla w|^2}{\partial\nu} \leq  c_1|\nabla w|^2~~\hbox{on} ~\partial\Omega~\hbox{for some}~ c_1>0
$$
and
$$
\||\nabla w|^2\|_{L^2(\partial\Omega)}\leq \eta\|\nabla|\nabla w|^2\|_{L^2(\Omega)}+c_2(\eta)\||\nabla w|^2\|_{L^1(\Omega)}~\hbox{for any}~ \eta\in(0,\frac54)
$$
(see Lemma 4.2 of \cite{Mizoguchi} and  Remark 52.9 in \cite{Quittner}), one can conclude that
\begin{equation}\label{3.18}
2 \displaystyle  \int_{\partial\Omega} |\nabla w|^2 \frac{ \partial|\nabla w|^2}{\partial\nu}\leq \frac 14\int_{\Omega}|\nabla |\nabla w|^2|^2
+c_3(\int_{\Omega}|\nabla w|^2)^2
\end{equation}
for some $c_3>0$.

For the other integrals on the right-side of \eqref{3.17}, we  use the Young inequality %and $ |\triangle w|^2\leq 2|D^2 w|^2$ on $\Omega$
to estimate
\begin{equation}\label{3.19}
-4 \int_{\Omega}|\nabla w|^2 \nabla w \cdot\nabla |\nabla w |^2
\leq
\displaystyle \frac{1}{3}\int_{\Omega}|\nabla |\nabla w|^2|^2+12 \int_{\Omega}|\nabla w|^6
\end{equation}
\begin{equation}\label{3.20}
-4\displaystyle  \int_{\Omega}n \nabla|\nabla w|^2  \cdot \nabla w
\leq
\displaystyle \frac{1}{3}\int_{\Omega}|\nabla |\nabla w|^2|^2+12 \int_{\Omega}n^2|\nabla w|^2
\end{equation}
as well as
 \begin{equation}\label{3.21}
\begin{array}{rl}
-4\displaystyle\int_{\Omega}n |\nabla w|^2   \triangle w
&
\leq
\displaystyle\frac16 \int_{\Omega} |\nabla w|^2 |\triangle w|^2+24\int_{\Omega}n^2|\nabla w|^2 \\[2mm]
&\leq \displaystyle \frac 13 \int_{\Omega} |\nabla w|^2 |D^2 w|^2+24\int_{\Omega}n^2|\nabla w|^2
\end{array}
\end{equation}
due to $ |\triangle w|^2\leq 2|D^2 w|^2$ on $\Omega$.

Substituting \eqref{3.18}--\eqref{3.21} into  \eqref{3.17}, we readily get
\begin{equation*}
\begin{array}{ll}
& \displaystyle \frac{d}{dt}\int_{\Omega}|\nabla w|^4+ \frac {13}{12}\int_{\Omega}|\nabla |\nabla w|^2|^2+
\frac {11}3\int_{\Omega}|\nabla w|^2  |D^2 w|^2
\\[3mm]
\leq &\displaystyle
12\int_{\Omega} |\nabla w|^6+36 \int_{\Omega}n^2|\nabla w|^2+c_3(\int_{\Omega}|\nabla w|^2)^2+4\int_{\Omega} |\nabla w|^4 |\nabla u|
\end{array}
\end{equation*}
and thus
\begin{equation}\label{3.22}
\begin{array}{ll}
&\displaystyle \frac{d}{dt}\int_{\Omega}|\nabla w|^4+ 2\int_{\Omega}|\nabla |\nabla w|^2|^2\\
\leq &
\displaystyle
12\int_{\Omega} |\nabla w|^6+36 \int_{\Omega}n^2|\nabla w|^2+c_3(\int_{\Omega}|\nabla w|^2)^2+4\int_{\Omega} |\nabla w|^4 |\nabla u|
\end{array}
\end{equation}
due to the fact $  |\nabla |\nabla w|^2|^2\leq 4 |\nabla w|^2  |D^2 w|^2 ~\hbox{on}~\Omega $.

Therefore combining \eqref{3.16} with \eqref{3.22} leads to
\begin{equation}\label{3.23}
\begin{array}{ll}
& \displaystyle \frac{d}{dt}(\int_{\Omega}n^2+
\int_{\Omega}|\nabla w|^4)+ 2\int_{\Omega}|\nabla |\nabla w|^2|^2+
\displaystyle\int_\Omega |\nabla n|^2
\\[3mm]
\leq &\displaystyle
12\int_{\Omega} |\nabla w|^6+37 \int_{\Omega}n^2|\nabla w|^2+c_3(\int_{\Omega}|\nabla w|^2)^2+2r \int_\Omega n^2- 2\mu \int_\Omega n^3
+4\int_{\Omega} |\nabla w|^4 |\nabla u|\\
\leq &\displaystyle
13\int_{\Omega} |\nabla w|^6+37^2 \int_{\Omega}n^3+c_3(\int_{\Omega}|\nabla w|^2)^2+2r \int_\Omega n^2- 2\mu \int_\Omega n^3+4\int_{\Omega} |\nabla w|^4 |\nabla u|.
\end{array}
\end{equation}
Furthermore by Lemma \ref{Lemma2.2} (ii), we get %from due to  the Gagliardo--Nirenberg inequality  in two-dimensional settings:
$ \|\varphi\|^3_{L^3}\leq K_2\|\nabla \varphi\|^2_{L^2}\|\varphi\|_{L^1}+ c_4\|\varphi\|^3_{L^1}$ and thus
$$
\int_{\Omega} |\nabla w|^6\leq K_2(\int_{\Omega}|\nabla |\nabla w|^2|^2)\left(\int_{\Omega} |\nabla w|^2\right)+c_4(\int_{\Omega} |\nabla w|^2)^3.
$$

Upon inserting this into \eqref{3.23} and \eqref{3.14}, we obtain
\begin{equation*}\label{3.24}
\begin{array}{ll}
& \displaystyle \frac{d}{dt}(\int_{\Omega}n^2+
\int_{\Omega}|\nabla w|^4)+ (2-13K_2\int_{\Omega} |\nabla w|^2)\int_{\Omega}|\nabla |\nabla w|^2|^2+
\displaystyle\int_\Omega |\nabla n|^2+  \int_{\Omega}n^2+
\int_{\Omega}|\nabla w|^4
\\[3mm]
\leq &\displaystyle
37^2 \int_{\Omega}n^3+(2r+1) \int_\Omega n^2- 2\mu \int_\Omega n^3+
\displaystyle \int_{\Omega}|\nabla w|^4+4\int_{\Omega} |\nabla w|^4 |\nabla u|+c_5,
\end{array}
\end{equation*}
which, along with
 $$
\int_{\Omega}|\nabla w|^4\leq \frac 17\int_{\Omega}|\nabla |\nabla w|^2|^2 +c_6
$$
and
\begin{equation*}
\begin{array}{ll}
4 \displaystyle\int_{\Omega} |\nabla w|^4 |\nabla u|
&\leq 4\||\nabla w|^2\|^2_{L^6(\Omega)}\|\nabla u\|_{L^{\frac 32}(\Omega)}\\
&\leq \displaystyle \frac{13}{56} \int_{\Omega}|\nabla |\nabla w|^2|^2 +c_7\\
\end{array}
\end{equation*}
by the Gagliardo--Nirenberg inequality and \eqref{2.13},
implies that
\begin{equation}\label{3.25}
\begin{array}{ll}
& \displaystyle \frac{d}{dt}(\int_{\Omega}n^2+
\int_{\Omega}|\nabla w|^4)+ (\frac {13}8-13K_2\int_{\Omega} |\nabla w|^2)\int_{\Omega}|\nabla |\nabla w|^2|^2+
\displaystyle\int_\Omega |\nabla n|^2+  \int_{\Omega}n^2+
\int_{\Omega}|\nabla w|^4
\\[3mm]
\leq &\displaystyle
37^2 \int_{\Omega}n^3+(2r+1) \int_\Omega n^2- 2\mu \int_\Omega n^3+c_8.
\end{array}
\end{equation}
Finally according to an extended variant (\cite{Biler}),
\eqref{3.6} and \eqref{3.14},
one can infer that
$$
\begin{array}{rl}
37^2 \displaystyle\int_{\Omega}n^3\leq &
 c_9\displaystyle\left(\int_\Omega |\nabla n|^2\right)\left(\int_\Omega  n |\ln n|\right)+c_9
(\int_\Omega  n)^3+c_9\\[2mm]
\leq &\displaystyle\frac 12 \int_\Omega |\nabla n|^2+c_{10}.
\end{array}
$$
Hence, in conjunction with \eqref{3.14}, \eqref{3.23} yields $c_{11}>0$ such that for all $t>t_*$
\begin{equation}\label{3.25}
 \displaystyle \frac{d}{dt}\int_{\Omega}(n^2+
|\nabla w|^4)+
\displaystyle  \int_{\Omega} (n^2+
|\nabla w|^4) +(\frac {13}8-13K_2\int_{\Omega} |\nabla w|^2)\int_{\Omega}|\nabla |\nabla w|^2|^2
\leq  c_{11}.
\end{equation}

Now in view of  \eqref{3.14}, we can see that
$\int_{\Omega}|\nabla w|^2\leq \frac r \mu(\frac{18}\chi+\frac 8{\chi^2} )|\Omega|$ for all $t>t_*$ and thus
$$
\frac {13}7-13K_2\int_{\Omega} |\nabla w|^2\geq 0
$$
provided that $\frac r{\mu}< \frac{\chi^2}{182 K_2|\Omega|}$, which guarantees that
 \begin{equation}\label{3.24}
 \displaystyle \frac{d}{dt}\int_{\Omega}(n^2+
|\nabla w|^4)+
\displaystyle  \int_{\Omega} (n^2+
|\nabla w|^4) \leq  c_{11}
\end{equation}for  all $t>t_*$ and thereby \eqref{3.15} is valid.

We are now ready to
prove Theorem 1.1 in the case of $r>0$.

{\it Proof of  Theorem 1.1 in the case of  $r>0$.}~ From the above lemmas, it follows that there exists $C>0$ such that
\begin{equation}\label{3.27}
\displaystyle  \|n(\cdot,t)\|_{L^2(\Omega)}+\|\nabla w(\cdot,t)\|_{L^4(\Omega)}+\|\nabla u(\cdot,t)\|_{L^2(\Omega)}\leq   C~~
\end{equation}
whenever $\mu>\mu_0(\chi, \Omega, r):=\max\{\mu_2(\chi, \Omega, r), \frac{182 K_2|\Omega|r}{\chi^2}\}$.
%\frac{ \chi}{K_1(36+16\chi)|\Omega|}\} $.
Thereupon, by the argument in e.g.\ Lemma 4.4 of \cite{Black},
%based  oadapting the  standard smoothing estimates for the Neumann heat semigroup,
we can readily prove that $\|n(\cdot,t)\|_{L^\infty(\Omega)},\|\nabla w(\cdot,t)\|_{L^\infty(\Omega)}$ and $\|A^\alpha u(\cdot,t)\|_{L^2(\Omega)}$ with some $\alpha\in (\frac 12,1)$ are globally bounded;
 %with respect  to the norm in $L^\infty(\Omega)$,
 we refer the reader to the proof of Lemma 4.4  in \cite{Black}, Lemma 3.12  and Lemma 3.11 in \cite{TaoZAMP} for the details. %, and thus omit it here.

Based on the global boundedness of solutions, we are able to  derive the convergence result claimed  in Theorem 1.1, namely,
\begin{equation}\label{3.28}
\displaystyle\lim_{t\rightarrow \infty}\|n(\cdot,t)-\frac r \mu\|_{L^\infty(\Omega)}=0,
\end{equation}
\begin{equation}\label{3.29}
 ~~\displaystyle\lim_{t\rightarrow \infty}\|\nabla w(\cdot,t)\|_{L^\infty(\Omega)}=0,
\end{equation}
\begin{equation}\label{3.29a}
 ~~\displaystyle\lim_{t\rightarrow \infty}\|u (\cdot,t)\|_{L^\infty(\Omega)}=0
\end{equation}
as well as
\begin{equation}\label{3.30}
 ~~\displaystyle\lim_{t\rightarrow \infty}\displaystyle\inf_{x\in \Omega } w(x,t) =\infty .
\end{equation}
 In fact, due to  $$
%\|n(\cdot,t)\|_{L^\infty(\Omega)}+
\int^{\infty}_{t_*}\int_\Omega  \frac{ |\nabla n|^2}n+
 \int^{\infty}_{t_*}\int_\Omega |\triangle w|^2\leq  C
 $$
 established in \eqref{3.14},  we can show \eqref{3.29}, \eqref{3.29a} and
\begin{equation}\label{4.4}
\displaystyle\lim_{t\rightarrow \infty}\|n(\cdot,t)-\overline n(t)\|_{L^\infty(\Omega)}=0
\end{equation}
with $\overline n(t)=\frac 1 {|\Omega|}\int_\Omega n(\cdot, t)$ by the arguments in Proposition 4.15 of \cite{Black}. Therefore it suffices to show that
\begin{equation}\label{3.33}
\displaystyle\lim_{t\rightarrow \infty}|\overline n(t)-\frac r\mu |=0.
\end{equation}
To this end, % similar to  that %inspired by we and give the skip of  for the conveentc.
we adapt the idea  of \cite{Litcanu} and  give the details of the proof for the convenience of readers.

Integrating the first equation  in \eqref{2.1} on the spatial variable over $\Omega$, we  obtain
$$
\overline n_t=  \lambda\overline n - \frac  \mu {|\Omega|} \int_\Omega n^2 =\lambda \overline n -\mu \overline n^2 - \frac {\mu}{|\Omega|}\int_\Omega (n-\overline n)^2.
$$
Putting $a(t):= \frac {\mu}{|\Omega|}\int_\Omega (n(\cdot,t)-\overline n)^2 $, the above equation then becomes
\begin{equation}\label{4.6a}
\overline n_t=  \mu \overline n(\frac \lambda \mu - \overline n) - a(t)
\end{equation}
Thereupon multiplying \eqref{4.6a} by $\overline n-\frac \lambda \mu$, we get
\begin{equation*}
\frac d{dt}(\overline n-\frac \lambda \mu)^2+ 2\mu \overline n(\overline n-\frac \lambda \mu )^2 =  -2 a(t)(\overline n-\frac \lambda \mu )
\end{equation*}
and then
\begin{equation}\label{3.36}
2\mu \int^\infty_{1}\overline n(\overline n-\frac \lambda \mu )^2\leq (\overline n(1)-\frac \lambda \mu)^2+
2 \sup_{t\geq 1}|\overline n(t)-\frac \lambda \mu | \int^\infty_{1}a(t).
\end{equation}
In addition, invoking the Poincar\'{e}--Wintinger inequality
$$\displaystyle\int_\Omega |\varphi-\frac 1 {|\Omega|}\int_\Omega \varphi(y)dy|^2  \leq C_p\int_\Omega |\varphi| \int_\Omega \frac{|\nabla \varphi|^2}{|\varphi|}  ~~\hbox{
for all~} \varphi\in W^{1,2}(\Omega)$$
for some $C_p > 0$, one can find
\begin{equation}\label{3.36a}
  \displaystyle\int^{\infty}_1 a(s)ds  \leq C_p\displaystyle\sup_{t\geq 1}\|n(t)\|_{L^1(\Omega)}
\displaystyle \int^{\infty}_1  \int_\Omega \frac{|\nabla n(s)|^2}{n(s)}  ds\leq C \\
\end{equation}
due to \eqref{3.14} and  Lemma 2.5(i). Hence combining \eqref{3.36a} with \eqref{3.36} yields
\begin{equation}\label{3.37}
\int^\infty_{1}\overline n(\overline n-\frac \lambda \mu )^2\leq C.
\end{equation}
On the other hand,
\begin{equation*}
\frac d{dt}\overline n(\overline n-\frac \lambda \mu )^2=\overline n_t((\overline n-\frac \lambda \mu)^2+2\overline n(\overline n-\frac \lambda \mu)),
\end{equation*}
which along with $|\overline n_t|\leq \lambda\overline n + \frac  \mu {|\Omega|} \int_\Omega n^2\leq C$ implies that
\begin{equation}\label{3.38}
\left|\frac d{dt}\overline n(\overline n-\frac \lambda \mu )^2\right|\leq C.
\end{equation}
Therefore by Lemma 6.3 of  \cite{Litcanu}, \eqref{3.38} and \eqref{3.37} show that
\begin{equation}\label{3.39}
\displaystyle\lim_{t\rightarrow \infty}   \overline n(t)(\overline n(t)-\frac r\mu)^2=0.
\end{equation}

From \eqref{4.4}, it follows that there exists $t_1>t_*$  such that
$
\|n(\cdot,t)-\overline n(t)\|_{L^\infty(\Omega)}\leq \frac \lambda {2\mu}
$
for all $t>t_1$, and thus
\begin{equation}\label{3.40}
\begin{array}{rl}
\overline n_t= & \lambda \displaystyle\overline n -\mu \overline n^2 - \displaystyle\frac {\mu}{|\Omega|}\int_\Omega n(n-\overline n)\\
\geq & \mu \overline n(\displaystyle\frac \lambda {\mu}-\overline n-\sup_{t>t_1}\|n(\cdot,t)-\overline n(t)\|_{L^\infty(\Omega)})\\
\geq & \mu \overline n(\displaystyle\frac \lambda {2\mu}-\overline n).
\end{array}
\end{equation}
On the other hand, noticing that the solution  $y(t)$ of the ODE
$$
y'(t)=\mu \overline y(\displaystyle\frac \lambda {2\mu}-\overline y), ~~y(t_1)>0
$$
satisfies
 $
\displaystyle\lim_{t\rightarrow \infty}  y(t)=\frac r {2\mu},
$
by the comparison principle, \eqref{3.40} implies that there exists $t_2>t_1$ such that for all $t\geq t_2$,
$$\overline n(t)\geq \displaystyle\frac \lambda {4\mu}.
$$
This together with \eqref{3.39} yields \eqref{3.33}.

Finally, in view of \eqref{3.28}, one can find  $t_3>1$ such that $n(x,t)\geq\frac r {2\mu}$
for all $x \in \Omega $  and $t\geq t_3$, and thereby $w(x,t)$ satisfies
$$w_t\geq \triangle w- |\nabla w|^2 +\frac r {2\mu}-u\cdot \nabla w $$
for $t\geq t_3$. Hence if $y(t)$ denotes the solution of ODE:
$y'(t)=\frac r {2\mu},~~y(t_3)=\displaystyle\min_{x\in \Omega}w(\cdot,t_3)$,
 then \begin{align}\label{4.13a}
 w(x,t)\geq \frac r {2\mu}(t-t_3)\end{align}
 by means of a straightforward parabolic comparison
which warrants that  \eqref{3.30} holds and thereby completes the proof.

\subsection{The case $r\leq 0$}
In this subsection, we show the global boundedness of solutions to \eqref{1.1}, \eqref{1.4}, \eqref{1.5} in the case $r\leq 0,\mu>0$. As mentioned in the introduction, due to the structure of \eqref{1.1} with  $r\leq 0,\mu>0$,  it is difficult to find a decreasing energy functional compared with the situation
when $r>0,\mu>0$ considered in the previous subsection or when $r=\mu=0$  considered in \cite{Winklerp}.  % is lacking a kind of structure
 %which  makes the availability  of the decreasing energy functional.
 % Let us note that the functional $\mathcal{F}(u,w)$ is invalid and more importantly,  the argument  in \cite{Winklerp} can not used directly in present situation.
Indeed, the  energy-type functional  $\mathcal{F}(n,w)$ in (3.1) of \cite{Winklerp} decreases along a solution in $\Omega\times (t_0,\infty)$ if $\mathcal{F}(n(\cdot,t_0),w(\cdot,t_0))$ is suitably small, namely
$$\frac d{dt}{\mathcal{F}(n,w)}\leq 0 ~~\hbox{for all}~~ t\geq t_0.
$$
%Accordingly, we make use of \eqref{2.2} to identify that the  functional of form
The main idea underlying our approach is to make use of the quadratic degradation in the first equation
of $\eqref{1.1}$ which should enforce some suitable regularity properties.  More precisely, on the basis of  \eqref{2.2},
 we can show  that the quantity of form
\begin{equation}\label{3.42}
\mathcal{F}(n,w):= \int_\Omega n(\ln n+a) dx+\frac \chi2\displaystyle \int_\Omega |\nabla w|^2 dx,
 \end{equation}
with parameter $a>0$ determined below (see \eqref{3.48}),  satisfies a certain of differential inequality. Although unlike the case of $r>0$
in which it  enjoys the monotonicity property, $\mathcal{F}(n,w)$ also provides us  the global boundedness of
$\int_\Omega n|\ln n| dx $ and $\int_\Omega |\nabla w|^2 dx$. This is encapsulated in the following lemma.

%The main ingredient is the analysis to energy-type inequalities.

\begin{lemma}\label{lemma34}
Let $\Omega\subset \mathbb{R}^2$ be a smooth bounded domain and $(n,w,u)$ be the global classical solution \eqref{2.1} with $r\leq 0, \mu>0$. Then there exists $t_*>0$
such that for all $t> t_*$
\begin{equation}
\label{3.43}
\int_\Omega |\nabla w(\cdot,t)|^2\leq \frac1{4 K_1} %\min\{\frac1{4 K_1}, \frac1{7 K_2}\}
\end{equation}
as well as
\begin{equation}\label{3.44}
\int_\Omega n|\ln n|  \leq C
\end{equation}
for some $C>0$.
\end{lemma}
\it Proof.  \rm  We test the first equation in \eqref{2.1} against $\ln n+a+1$, and  integrate by parts to see that
\begin{align}
&\displaystyle\frac{d}{dt}\int_\Omega n(\ln n+a)\nonumber\\
\leq & -\displaystyle \displaystyle\int_\Omega  \frac{ |\nabla n|^2}n-\chi\displaystyle\int_\Omega \nabla n\cdot\nabla w +
 \displaystyle\displaystyle\int_\Omega (n(r- \mu n)-u\cdot \nabla n) (\ln n+a+1)
 \nonumber\\
 \leq & -\displaystyle \displaystyle\int_\Omega  \frac{ |\nabla n|^2}n-\chi\displaystyle\int_\Omega \nabla n\cdot\nabla w +
 \displaystyle\displaystyle\int_\Omega n(r- \mu n) (\ln n+a).\label{3.45}
\end{align}
 due to $r\leq 0$ and $\nabla\cdot u=0$.

 On the other hand, recalling  \eqref{3.4} and \eqref{2.14}, it is possible to fix $t_0>0$ such that for all $t\geq t_0$, we have
 \begin{equation}\label{3.46}
 \begin{array}{ll}
&\displaystyle\frac 1 2\frac{d}{dt} \int_\Omega |\nabla w|^2+
\displaystyle\frac 1 4\int_\Omega |\triangle w|^2+
\displaystyle\frac 1 4(\frac 34-2K_1\|\nabla w\|^2_{L^2(\Omega)} ) \int_\Omega |\triangle w|^2
\\
\leq &  \displaystyle\int_\Omega \nabla u\cdot\nabla w.
\end{array}
\end{equation}
From Lemma \ref{Lemma2.2}(i),  %$\frac{\partial w}{\partial \nu}|_{\partial\Omega}=0 $ and the two-dimensional
%version of the Gagliardo--Nirenberg inequality,
there exists a constant $K_3>0$ such that
\begin{align}
\label{3.47a} 8 K_3\|\nabla w \|^2_{L^2(\Omega)}\leq\|\Delta w\|^2_{L^2(\Omega)}.
\end{align}
%which along with \eqref{3.31} yields
 %\begin{equation}\label{3.32}
% \displaystyle\frac 1 2\frac{d}{dt} \int_\Omega |\nabla w|^2+\displaystyle\frac 1 4\displaystyle \int_\Omega |\triangle w|^2+ K_3
%(1-2 K_1\|\nabla w\|^2_{L^2(\Omega)} ) \int_\Omega |\nabla w|^2
%\leq
% \displaystyle\int_\Omega \nabla u\cdot\nabla w.
%\end{equation}
Hence combining \eqref{3.46} with \eqref{3.45}, we get
\begin{align}
&\displaystyle\frac d{dt}{\mathcal{F}(n,w)} +\displaystyle \displaystyle\int_\Omega  \frac{ |\nabla n|^2}n+\displaystyle\frac \chi 4\int_\Omega |\triangle w|^2+ K_3
\int_\Omega n(\ln n+a)\nonumber\\
&+\frac \chi 4
(\frac 3 4-2 K_1\|\nabla w\|^2_{L^2(\Omega)} ) \int_\Omega |\triangle w|^2\nonumber\\
\leq
&
\displaystyle\int_\Omega n(K_3- \mu n)(\ln n+a)+r\displaystyle\int_\Omega n(\ln n+a)~~\hbox{for }~t\geq t_0.\label{3.47}
\end{align}
Now for any fixed $\varepsilon<\min \{\frac  \chi{24K_1}, \frac  \chi{42K_2}\}$, we pick $a>1$ sufficiently large  such that
\begin{align}\label{3.48}
e^{-a}< \frac {K_3}\mu, ~~~(1-r)|\Omega|\displaystyle\max_{0<n\leq e^{-a}} |n\ln n|< \varepsilon \displaystyle\min\{K_3,1\},
\end{align}
due to  $\displaystyle\lim_{n \rightarrow 0} n\ln n=0$ and $n\ln n<0 $
for all $n\in (0,1)$, and thereby fix $t_1>\max\{1,t_0\}$ fulfilling
\begin{align}\label{3.49}
~~\displaystyle\frac {a|\Omega|}{\mu(t_1+\gamma)}<\frac\varepsilon {4},
\displaystyle\frac{|\Omega|}{\mu^2(t_1+\gamma)}<\frac\varepsilon {16},
\displaystyle\frac{(a+(\ln \frac{K_3} \mu)_+ )|\Omega|}{\mu (t_1+\gamma)}+\frac{2\chi(
\frac{|\Omega|}{\mu}+\|n_0\|_{L^1(\Omega)}+\|w_0\|_{L^1(\Omega)})}{t_1}< \displaystyle\frac \varepsilon 4.
\end{align}

Let $t_2=t_1+t_1^2$,
 $$\mathcal{S}_1\triangleq\{t\in[t_1,t_2]| \int_\Omega |\nabla w(\cdot,t)|^2\geq\frac \varepsilon{2\chi} \}
 $$
 and
$$\mathcal{S}_2\triangleq\{t\in[t_1,t_2]| \|n(\cdot,t)\|^2_{L^2(\Omega)}\geq  \frac \varepsilon 4 \}.
$$
Then \begin{equation}\label{3.50}
 |\mathcal{S}_1|\leq \frac {|t_2-t_1|}4, ~~~|\mathcal{S}_2|\leq \frac {|t_2-t_1|}4.
 \end{equation}
 %To verify the claim,
  By Lemma \ref{lemma23}(iii), \eqref{2.2} and the second equation in \eqref{2.1}, we obtain that
\begin{align}
\displaystyle  \int^{t_2}_{t_1}  \int_\Omega |\nabla w|^2 \leq &
\int^{t_2}_{t_1}  \int_\Omega n  +
 \int_\Omega w(\cdot,t_1)\nonumber \\
=& \int^{t_2}_{t_1}  \int_\Omega n  +\int_\Omega w_0+\int^{t_1}_{0}  \int_\Omega n\nonumber\\
\leq & \displaystyle\frac {|\Omega|}{\mu(t_1+\gamma)}(t_2-t_1)+ \int_\Omega w_0+t_1\int_\Omega n_0\nonumber.
\end{align}
Furthermore, by \eqref{3.49}
\begin{align*}
\displaystyle  \int^{t_2}_{t_1}  \int_\Omega |\nabla w|^2 dxds\leq &
\displaystyle(\frac {|\Omega|}{\mu(t_1+\gamma)}+ \frac{t_1\|n_0\|_{L^1(\Omega)}+\|w_0\|_{L^1(\Omega)}}{t_2-t_1})(t_2-t_1)\nonumber\\
\leq &
\frac{\frac{|\Omega|}{\mu}+\|n_0\|_{L^1(\Omega)}+\|w_0\|_{L^1(\Omega)}}
{t_1}(t_2-t_1)\\
< & \frac{\varepsilon}{8\chi}(t_2-t_1).
\end{align*}
On the other hand, by the definition of $\mathcal{S}_1$, we see that
$$\frac \varepsilon{2\chi}|\mathcal{S}_1|\leq
\displaystyle  \int^{t_2}_{t_1}  \int_\Omega |\nabla w|^2.
$$
 and thereby  $|\mathcal{S}_1|\leq \frac {|t_2-t_1|}4$.

 In addition, by \eqref{2.2} and \eqref{3.49}, we get
$$\displaystyle  \int^{t_2}_{t_1}  \int_\Omega n^2 \leq \frac 1\mu   \int_\Omega n(\cdot,t_1)\leq \displaystyle\frac{|\Omega|}{\mu^2(t_1+\gamma)}<\frac\varepsilon {16},$$
which implies that $|\mathcal{S}_2|\leq \frac {|t_2-t_1|}4$.

Therefore from \eqref{3.50}, it follows that
$
|(t_1,t_2)\setminus (\mathcal{S}_1\cup  \mathcal{S}_2)| \geq  \frac {|t_2-t_1|}2,
$ and thereby  there exists $t_*\in (t_1,t_2) $  such that
\begin{equation}\label{3.51}
\displaystyle \|n(\cdot,t_*)\|^2_{L^2(\Omega)}<  \frac\varepsilon 4
\end{equation}
and
\begin{equation}\label{3.52}
\int_\Omega |\nabla w(\cdot,t_*)|^2<\frac {\varepsilon} {2\chi}<\frac {1} {6K_1}.
\end{equation}
By \eqref{3.52}, we can see that the set
$$
\mathbf{S}\triangleq \{ t\in ( t_*,\infty)| ~K_1\int_\Omega |\nabla w(\cdot,t)|^2<\frac14 \}
$$
is not empty and hence $T_S=\sup {\mathbf{S}}$ is a well-defined element of  $(t_*,\infty]$. In fact, we claim that $T_S=\infty$.
To this end, supposing on the contrary that  $T_S<\infty$, we then have $K_1\int_\Omega |\nabla w(\cdot,t)|^2<\frac14$ for all $t\in [t_*,T_S)$, but
\begin{equation}\label{3.54}
K_1\int_\Omega |\nabla w(\cdot,T_S)|^2=\frac14.
\end{equation}
Hence from %and, by the definition of $T_S$ and
\eqref{3.47} and  \eqref{3.47a}, it follows that for all $t\in [t_*,T_S)$,
 \begin{align}
&\displaystyle\frac d{dt}{\mathcal{F}(n,w)} +\displaystyle \displaystyle\int_\Omega  \frac{ |\nabla n|^2}n+\displaystyle\frac \chi 4\int_\Omega |\triangle w|^2+K_3
\int_\Omega n(\ln n+a)+\frac {K_3\chi}2
\int_\Omega |\nabla w|^2\nonumber\\
\leq
&
\displaystyle\int_\Omega n(K_3- \mu n)(\ln n+a)+ r \int_\Omega n(\ln n+a)\nonumber\\
\leq
&
\displaystyle\int_{e^{-a}<n\leq \frac {K_3} \mu}n(K_3- \mu n)(\ln n+a)+
\displaystyle r\int_{0<n\leq e^{-a}}n(\ln n+a)\nonumber\\
\leq
&
K_3\displaystyle\int_{e^{-a}<n\leq \frac{K_3} \mu}n(\ln n+a)+ \displaystyle r\int_{0<n\leq e^{-a}}n\ln n \label{3.55}\\
\leq
&
aK_3\displaystyle\int_\Omega n+K_3
\displaystyle\int_{e^{-a}<n\leq \frac{K_3} \mu}n\ln n
-r|\Omega|\displaystyle\max_{0<n\leq e^{-a}} |n\ln n| \nonumber\\
\leq
&K_3
(a+(\ln \frac{K_3} \mu)_+ )\displaystyle\int_\Omega n + \varepsilon K_3\nonumber\\
\leq
&
\displaystyle\frac{(a+(\ln \frac{K_3} \mu)_+ )K_3|\Omega|}{\mu(t_1+\gamma)}+ \varepsilon K_3,\nonumber
\end{align}
where  we have made use of  $t_*\geq t_1$, the decay estimate \eqref{2.2} and \eqref{3.49},  and thus
 \begin{align*}
 & \mathcal{F}(n,w)(T_s)+\displaystyle\int^{T_s}_{t_*} \displaystyle e^{- K_3(T_s-\sigma)}(\int_\Omega \frac{ |\nabla n|^2}n(\cdot,\sigma)+\displaystyle\frac \chi 4\int_\Omega |\triangle w(\cdot,\sigma)|^2)d\sigma\\
 \leq &  \mathcal{F}(n,w)(t_*)+\displaystyle\frac{(a+(\ln \frac{K_3} \mu)_+ )|\Omega|}{\mu(t_1+\gamma)}+ \varepsilon,
\end{align*}
which implies that
 \begin{align}
\displaystyle\frac \chi 2\int_\Omega |\nabla w(\cdot,T_S)|^2 \leq &
\mathcal{F}(n,w)(t_*)+\displaystyle\frac{(a+(\ln \frac{K_3} \mu)_+ )|\Omega|}{\mu (t_1+\gamma)}-\int_\Omega n(\ln n+a)(\cdot,T_S)+ \varepsilon\nonumber\\
\leq & \int_\Omega n(\ln n+a) (\cdot,t_*)+\frac \chi2\displaystyle \int_\Omega |\nabla w|^2(\cdot,t_*)+ \varepsilon\nonumber\\
&+\displaystyle\frac{(a+(\ln \frac{K_3} \mu)_+ )|\Omega|}{\mu (t_1+\gamma)}-\int_\Omega n(\ln n+a)(\cdot,T_S)\nonumber\\
\leq & \int_\Omega (n^2+an) (\cdot,t_*)+\frac \chi2\displaystyle \int_\Omega |\nabla w|^2(\cdot,t_*)+ \varepsilon\label{3.56}\\
&+\displaystyle\frac{(a+(\ln \frac{K_3} \mu)_+ )|\Omega|}{\mu (t_1+\gamma)}-\int_\Omega n(\ln n+a)(\cdot,T_S),\nonumber
\end{align}
due to $n\geq \ln n $ for all $n>0$.

In addition, by \eqref{3.49}, we see that
\begin{align}
\int_\Omega n(\ln n+a)(\cdot,T_S)\geq &\int_{0<n\leq e^{-a}} n(\ln n+a)(\cdot,T_S)\label{3.57}\\
\geq &\int_{0<n\leq e^{-a}} n\ln n (\cdot,T_S)\nonumber\\
\geq & -|\Omega|\displaystyle\max_{0<n\leq e^{-a}} |n\ln n|\nonumber\\
\geq & -\varepsilon \nonumber.
\end{align}
Upon inserting \eqref{3.57} into \eqref{3.56}, we see that
 \begin{align}
\displaystyle\frac \chi 2\int_\Omega |\nabla w(\cdot,T_S)|^2 \leq &
\int_\Omega (n^2+an) (\cdot,t_*)+\frac \chi2\displaystyle \int_\Omega |\nabla w|^2(\cdot,t_*)\label{3.58}\\
&+\displaystyle\frac{(a+(\ln \frac{K_3} \mu)_+ )|\Omega|}{\mu(t_1+\gamma)}+2\varepsilon,\nonumber
\end{align}
which along with \eqref{3.51}, \eqref{3.52}, \eqref{2.2} and  \eqref{3.49}, establishes that
\begin{align}\label{3.59}
\displaystyle\frac \chi 2\int_\Omega |\nabla w(\cdot,T_S)|^2 \leq &
\displaystyle \frac {5\varepsilon} 2+
a\int_\Omega n (\cdot,t_*)+\displaystyle\frac{(a+(\ln \frac{K_3} \mu)_+ )|\Omega|}{\mu (t_1+\gamma)}\nonumber\\
< &
3\varepsilon\\
\leq & \frac {\chi}{8K_1}.
\nonumber\end{align}
This contradicts  \eqref{3.54} and thereby $T_S=\infty$, which means that the differential inequality \eqref{3.55} is actually valid for all $t>t_*$.

Now revisiting the proof of \eqref{3.59}, upon integration in time over $(t_*,t)$, we  have
$$
\displaystyle\frac \chi 2\int_\Omega |\nabla w(\cdot,t)|^2 \leq 3\varepsilon~~~\hbox{for all}~~t>t_*
$$ which implies that  \eqref{3.43} is valid by the choice of $\varepsilon$, as well as
\begin{equation}\label{3.60}
\int_\Omega n \ln n(\cdot,t)  \leq C_1~~~\hbox{for all}~~t>t_*
\end{equation}
for some $C_1>0$.

Since $\xi\ln \xi\geq -\frac 1e $ for all $\xi>0$,
 \begin{align*}\label{3.45}
\int_\Omega n |\ln n|(\cdot,t) &=  \int_\Omega n \ln n(\cdot,t)-2 \int_{0<n<1} n \ln n(\cdot,t)\\
& \leq \int_\Omega n \ln n(\cdot,t)+\frac{2|\Omega|}{e},
\end{align*}
which along with \eqref{3.60}
 readily  implies that \eqref{3.44} is actually  valid with $C=C_1+\frac{2|\Omega|}{e}$.

Furthermore, from \eqref{3.55}, one can  also conclude that:

\begin{corollary}\label{Cor32}  %Let $(p,c,w)$ be the global classical solution of (1.1). Then
Under the conditions of Lemma \ref{lemma34}, we have
\begin{equation}\label{3.61}
\displaystyle\lim_{t\rightarrow \infty}\int^{t+1}_{t}\int_\Omega ( \frac{ |\nabla n|^2}n+
|\triangle w|^2) =0, \quad
\displaystyle\lim_{t\rightarrow \infty}\int_\Omega |\nabla w(\cdot,t)|^2=0.
\end{equation}
\end{corollary}
\it Proof.  \rm  On the basis of the decay estimate \eqref{2.2} and revisiting the argument in the proof of Lemma \ref{lemma34},  one can conclude  that
for any fixed  $\varepsilon\in (0,\frac  \chi{42K_1})$, there exists  $t_\varepsilon>1$ such that
%$\int_\Omega |\nabla w(\cdot,t)|^2<\varepsilon$ whenever $t>t_\varepsilon$. On the other hand,
%\begin{align*}
% &
$$
\int_\Omega |\nabla w(\cdot,t)|^2+\displaystyle\int^{t}_{t_\varepsilon} \displaystyle e^{- K_3(t-\sigma)}(\int_\Omega \frac{ |\nabla n|^2}n(\cdot,\sigma)+\displaystyle\frac \chi 8\int_\Omega |\triangle w(\cdot,\sigma)|^2)d\sigma
%\\
 \leq
%&
 \varepsilon
%\end{align*}
$$
for all $t>t_\varepsilon$.
Furthermore, it follows from the above inequality that
%\begin{align*}
% &
$$
\displaystyle\int^{t}_{t-1} \displaystyle(\int_\Omega \frac{ |\nabla n|^2}n(\cdot,\sigma)+\displaystyle\frac \chi 8\int_\Omega |\triangle w(\cdot,\sigma)|^2)d\sigma
%\\
 \leq
%&
\varepsilon e^{ K_3}
%\end{align*}
$$
for any $t>t_\varepsilon+1$, which implies that \eqref{3.61} is indeed valid.

At this point, we can
prove Theorem 1.1 in the case of $r\leq 0$.

{\it Proof of  Theorem 1.1 in the case $r\leq 0$.}~  We can repeat the argument in the proof of Theorem 1.1 in the case $r>0$.  In fact, in view of \eqref{3.43} and \eqref{3.44}, \eqref{3.15} is also valid for $r\leq 0,\mu>0$, and thereby the global boundedness of solutions can be  proven. In addition,   similar to  the case of $r>0$,   we  can show % as well as the convergence properties claimed in Theorem 1.1
\begin{equation}\label{3.62}
\displaystyle\lim_{t\rightarrow \infty}\|n(\cdot,t)\|_{L^\infty(\Omega)}=0,
\end{equation}
\begin{equation}\label{3.63}
 ~~\displaystyle\lim_{t\rightarrow \infty}\|\nabla w(\cdot,t)\|_{L^\infty(\Omega)}=0
 \end{equation}
 as well as
\begin{equation}\label{3.64}
\displaystyle\lim_{t\rightarrow \infty}\|u(\cdot,t)\|_{L^\infty(\Omega)}=0.
\end{equation}
%Based on \eqref{2.2} and \eqref{3.46}, \eqref{3.60} and \eqref{3.61}ed
%by the argument similar employed in the proof of Lemma 6.1 of \cite{Winklerp}.
For  the sake of completeness we shall only recount the main steps and refer to the mentioned sources for more details. Invoking standard parabolic regularity theory (see the proofs of Lemma 4.5 and Lemma 4.9 of \cite{Winklerp} for details), one can see that there exist $\theta\in(0,1)$ and $\alpha\in (\frac 12,1)$  and $c_1>0$ such that
for all $t>1$
\begin{equation}\label{3.65}
\|n\|_{C^{\theta,\frac \theta 2}(\overline{\Omega}\times [t,t+1])}+\|\nabla w(\cdot,t)\|_{C^{\theta}(\overline{\Omega})}+
 \|A^\alpha u(\cdot,t)\|_{L^2(\Omega)}
\leq c_1.
\end{equation}
If \eqref{3.62} were false, then there would be $c_2>0$, $(t_k)_{k\in\mathbb{N}}$  and  $(x_k)_{k\in\mathbb{N}}\subseteq \Omega$  such that $t_k\rightarrow \infty$ as $k\rightarrow \infty$, and $n(x_k,t_k)>c_2$ for all $k\in \mathbb{N}$, which, along with the uniform continuity of $n$ in $\overline{\Omega}\times [t,t+1]$  as shown by \eqref{3.65}, entails that one can find $r>0$ such that $B(x_k, r)\subseteq \Omega$ for all $k\in\mathbb{N}$ and
$$n(x,t_k)>\frac {c_2}2~~\hbox{for all}~~x\in B(x_k, r).$$
This shows $$
\int_{\Omega}n(\cdot,t_k)\geq \int_{B(x_k, r)}n(\cdot,t_k)\geq  \frac {c_2}2 \pi r^2
$$
 which contradicts \eqref{2.2} and thus proves \eqref{3.62}.  Similarly, on the basis of \eqref{3.61} and \eqref{3.65},   \eqref {3.63} can be proved.
 Finally, \eqref {3.64} results from \eqref{2.14}, \eqref{3.65} and a simple interpolation,  and  thereby completes the proof.

\section{Asymptotic profile of  solutions }
It is observed that in the case $r<0$,  solutions to \eqref{1.1}, \eqref{1.4}, \eqref{1.5} enjoy the exponential decay property due to the
exponential decay of $\| n(\cdot,t)\|_{L^1(\Omega)}$. Therefore the present paper focuses on the asymptotic profile of  \eqref{1.1}, \eqref{1.4}, \eqref{1.5} in the cases $r>0$ and $r=0$, namely, we will give the proofs of Theorems 1.2 and 1.3 respectively.

\subsection{The case $r>0$}  %Convergence of solutions with exponential rate }%Precisely,

Making use of the convergence properties of  $(n,\frac{|\nabla c|}c)$  asserted in Theorem 1.1, we apply
$L^p-L^q$ estimates  for the Neumann heat semigroup $(e^{t\Delta})_{t>0}$  to show
 %by a self-map-type reasoning
 %that whenever $\frac r \mu$ is suitably small,
  %combining with  the variation-of-constants representation of $u-\frac r \mu$ and $w$,  respectively,
%this information allows us to
%$$(u-\frac r \mu,v,\frac{|\nabla v|}v)\longrightarrow (0,0,0)$$
%On the basis of  the estimates of $(u,w)$ in \eqref{3.13} including $\int^\infty_{t_*}\int_\Omega  \frac{ |\nabla u|^2}u$, we shall  be able to show that
$(n,c,u)\rightarrow (\frac r \mu,0,0) $ in $L^\infty(\Omega)$ and $\frac {|\nabla c|}c \rightarrow 0$ in $L^p(\Omega)$ at some exponential rate as $t\rightarrow \infty$, respectively, whenever $\mu$ is suitably large compared with $r$.
To this end, we first make an observation which will be used in the proof of the subsequent lemma:

\begin{lemma}\label{Lemma41} For any $\alpha\in(0,\min\{\lambda_1, r\})$,  for $I:=\int_0^\infty(1+\sigma^{-\frac23}+\sigma^{-\frac12})e^{-(\lambda_1-\alpha) \sigma}d\sigma$,
and $c_i>0$ (i=1,3) as given by  Lemma 2.1,  there exist $\varepsilon_1>0,\varepsilon_2>0$ such that
\begin{align}
&4 c_1I \varepsilon_2<1, ~~4 c_1|\Omega|^{\frac 16} I \varepsilon_1\leq \varepsilon_2, ~~8\chi c_3 I\varepsilon_2<1,\label{4.1}\\
& 8 \mu \varepsilon_1<r-\alpha.  \label{4.2}
\end{align}
\end{lemma}
\begin{lemma}\label{lemma43} Let $(n,w,u)$ be the global bounded solution of \eqref{2.1}.  If  $\mu>  32 \chi c_3c_1|\Omega|^{\frac 16}I^2 r$,
one can find constants $C_i>0$, $i=1,2,3$, $\alpha\in(0,\min\{\lambda_1, r\})$ and $\beta<\alpha $ such that
\begin{align}\label{4.3}
\|n(\cdot,t)-\frac  r \mu\|_{L^\infty(\Omega)}&\leq C_1e^{-\alpha t},
\\
\|\nabla w(\cdot,t)\|_{L^6(\Omega)}&\leq C_2e^{-\alpha t}
\end{align}
as well as
\begin{align}\label{4.5}
\|u(\cdot,t)\|_{L^\infty(\Omega)}&\leq C_3 e^{-\beta t}
\end{align}
for all $t\geq 1 $.
% Here, $\lambda_1$ and $\lambda_1'$ are as defined in Lemmas 2.1 and 2.2 below.
\end{lemma}

\it{Proof.}\rm\quad Let $N(x,t)=n(x,t)-\frac  r \mu$, $\varepsilon_1>0 $ and $\varepsilon_2>0$ be given by Lemma \ref{Lemma41}.  Then from
\eqref{3.28}, \eqref{3.29} and \eqref{3.29a}, there exists $t_0>1$ suitably large such that
for $t\geq t_0$
\begin{align}\label{4.6}
\|N(\cdot,t)\|_{L^\infty(\Omega)}\leq \frac {\varepsilon_1} 8,~~~(c_2+1)\|\nabla w(\cdot,t)\|_{L^\infty(\Omega)} \leq \frac {\varepsilon_2} 8
\end{align}
and
\begin{align}\label{4.7}
8c_1\| u(\cdot,t)\|_{L^\infty(\Omega)}\int_0^\infty(1+\sigma^{-\frac12})e^{-(\lambda_1-\alpha) \sigma}d\sigma \leq 1.
\end{align}
Now  we consider
\begin{align}
T\!\triangleq\!\sup\!\left\{\!\widetilde{T}\!\in\!(t_0,\infty)\!\left|
\begin{array}{ll}
\!\| N(\cdot,t)\|_{L^\infty(\Omega)}\leq\varepsilon_1 e^{-\alpha (t-t_0)} \quad \hbox{for all }\, t\in[t_0,\widetilde{T}),
\\
\!\|\nabla w(\cdot,t)\|_{L^6(\Omega)}\leq\varepsilon_2 e^{-\alpha (t-t_0)} \quad \hbox{for all }\, t\in[t_0,\widetilde{T}).
\end{array}\right.
\!\right\}\label{4.8}
\end{align}
By \eqref{4.6}, $T$ is well-defined. In what follows, we shall demonstrate that $T=\infty$.

To this end,
we first invoke the  variation-of-constants representation of $w$:
\begin{equation}\label{4.9}
\begin{array}{ll}
 w(\cdot,t)=
 & e^{(t-t_0) \Delta }w(\cdot,t_0) -\displaystyle\int_{t_0}^t e^{(t-s)\Delta}|\nabla w(\cdot,s)|^2 ds+
\int_{t_0}^t e^{(t-s)\Delta}N(\cdot,s) ds\\[2mm]
 & \displaystyle-\int_{t_0}^t e^{(t-s)\Delta}(u\cdot \nabla w)(\cdot,s) ds+\frac r \mu (t-t_0),
 \end{array}
\end{equation}
and  use Lemma 2.1(i), (ii) to estimate
\begin{align}
& \|\nabla w(\cdot,t)\|_{L^6(\Omega)}\nonumber\\
\leq &
\displaystyle \| \nabla e^{(t-t_0) \Delta }w(\cdot,t_0)\|_{L^6(\Omega)}
+\displaystyle \int_{t_0}^t
\|
\nabla e^{(t-s)\Delta}|\nabla w(\cdot,s)|^2\|_{L^6(\Omega)} ds\nonumber\\
& +\displaystyle
\int_{t_0}^t \|\nabla e^{(t-s)\Delta}N(\cdot,s)\| _{L^6(\Omega)} ds
+\displaystyle\int_{t_0}^t \|\nabla e^{(t-s)\Delta}(u\cdot \nabla w)(\cdot,s)\|_{L^6(\Omega)}   ds\nonumber\\
 \leq & 2 c_2e^{-\lambda_1(t-t_0)}\|\nabla w(\cdot,t_0)\|_{L^6(\Omega)}+
c_1\displaystyle\int_{t_0}^t(1+(t-s)^{-\frac23})e^{-\lambda_1(t-s)}\|\nabla w(\cdot,s)\|^2_{L^6(\Omega)}ds\nonumber
\\
&+c_1|\Omega|^{\frac 16}\displaystyle\int_{t_0}^t(1+(t-s)^{-\frac12})e^{-\lambda_1(t-s)}\|N(\cdot,s)\|_{L^{\infty}(\Omega)}ds\nonumber\\
&
+
c_1\displaystyle\int_{t_0}^t(1+(t-s)^{-\frac12})e^{-\lambda_1(t-s)}\|u(\cdot,s)\|_{L^{\infty}(\Omega)}\|\nabla w(\cdot,s)\|_{L^6(\Omega)}ds\nonumber\\
:= &I_1+I_2+I_3\label{4.10}
\end{align}
for all $t_0<t<T$.

Now we estimate the terms $I_i$ $(i=1,2,3)$ respectively.  Firstly, from \eqref{4.6},
we have $I_1\leq \frac {\varepsilon_2} 4e^{-\lambda_1(t-t_0)}$. By the definition of $T$ and \eqref{4.1}, we can see that
\begin{align*}
\begin{array}{rl}
I_2\leq
&c_1\varepsilon^2_2\displaystyle\int_{t_0}^t(1+(t-s)^{-\frac23})e^{-\lambda_1(t-s)}e^{-2\alpha (s-t_0)}ds\\[2mm]
\leq &
c_1\varepsilon^2_2\displaystyle\int_{t_0}^t(1+(t-s)^{-\frac23})e^{-\lambda_1(t-s)}e^{-\alpha (s-t_0)}ds
\\[2mm]
\leq &
c_1\varepsilon^2_2\displaystyle\int_{0}^\infty (1+\sigma^{-\frac23})e^{-(\lambda_1-\alpha)\sigma}d \sigma\cdot e^{-\alpha (t-t_0)}\\
\leq &
 \displaystyle\frac {\varepsilon_2} 4e^{-\alpha(t-t_0)}.
\end{array}
\end{align*}
By the definition of $T$, \eqref{4.7} and \eqref{4.1} again, we also have
\begin{align*}
\begin{array}{rl}
I_3\leq
&(c_1|\Omega|^{\frac 16}\varepsilon_1+c_1\displaystyle\sup_{t\geq t_0}\| u(\cdot,t)\|_{L^\infty(\Omega)} \varepsilon_2)\displaystyle\int_{t_0}^t(1+(t-s)^{-\frac12})e^{-\lambda_1(t-s)}e^{-\alpha (s-t_0)}ds\\[2mm]
= &
(c_1|\Omega|^{\frac 16}\varepsilon_1+c_1\displaystyle\sup_{t\geq t_0}\| u(\cdot,t)\|_{L^\infty(\Omega)} \varepsilon_2)\displaystyle\int_{t_0}^t(1+(t-s)^{-\frac12})e^{-(\lambda_1-\alpha)(t-s)}e^{-\alpha (t-t_0)}ds
\\[2mm]
\leq &
(c_1|\Omega|^{\frac 16}\varepsilon_1+c_1\displaystyle\sup_{t\geq t_0}\| u(\cdot,t)\|_{L^\infty(\Omega)} \varepsilon_2)\displaystyle\int_{0}^\infty (1+\sigma^{-\frac12})e^{-(\lambda_1-\alpha)\sigma}d \sigma\cdot e^{-\alpha (t-t_0)}\\
\leq &
 \displaystyle\frac {3\varepsilon_2} 8e^{-\alpha(t-t_0)}.
\end{array}
\end{align*}
Substituting these estimates into \eqref{4.10}, we  get
\begin{align}\label{4.11}
\|\nabla w(\cdot,t)\|_{L^6(\Omega)}\leq \frac {7\varepsilon_2} 8 e^{-\alpha(t-t_0)}< \varepsilon_2 e^{-\alpha (t-t_0)} \quad \hbox{for all }\, t\in[t_0,T).
\end{align}

On the other hand, since $$
N_t=\triangle N+\chi \nabla\cdotp ( n\nabla w)-rN- \mu N^2- u\cdot\nabla N,
$$
 the variation-of-constants representation of $N$ yields
 $$
 \begin{array}{ll}
N(\cdot,t)= &\displaystyle e^{(t-t_0) (\Delta-r) }N(\cdot,t_0) +\chi \int_{t_0}^t e^{(t-s)(\Delta-r)}\nabla\cdotp ( n\nabla w)(\cdot,s)ds-
\mu\int_{t_0}^t e^{(t-s)(\Delta-r)}N^2(\cdot,s) ds
\\
&-\displaystyle\int_{t_0}^t e^{(t-s)(\Delta-r)}( u\cdot \nabla N)(\cdot,s)ds.
\end{array}
$$
 Then by $\nabla \cdot u=0$ we can see that
 \begin{align*}
\begin{array}{rl}
 & \|N(\cdot,t)\|_{L^{\infty}(\Omega)}\\
 \leq
 &
 \|\displaystyle e^{(t-t_0) (\Delta-r) }N(\cdot,t_0)\|_{L^{\infty}(\Omega)} +\displaystyle
\mu\int_{t_0}^t \| e^{(t-s)(\Delta-r)}N^2(\cdot,s) \|_{L^{\infty}(\Omega)}ds
 \\
 &+\displaystyle\int_{t_0}^t \| e^{(t-s)(\Delta-r)}\nabla\cdot ( u N)(\cdot,s)\|_{L^{\infty}(\Omega)}ds
  +\displaystyle\chi \int_{t_0}^t \| e^{(t-s)(\Delta-r)}\nabla\cdot ( n\nabla w)(\cdot,s)\|_{L^{\infty}(\Omega)}ds
\\
:=& J_1+J_2+J_3+J_4.\\
 \end{array}
\end{align*}
Here the maximum principle together with \eqref{4.6}  ensures  that
$$
J_1\leq e^{-r(t-t_0)} \|N(\cdot,t_0)\|_{L^\infty(\Omega)}\leq \frac {\varepsilon_1} 8 e^{-\alpha(t-t_0)}.
$$
By the definition of $T$ and  comparison principle,  we infer that
\begin{align*}
\begin{array}{rl}
J_2\leq &\displaystyle
\mu\int_{t_0}^te^{-r(t-s)} \|e^{(t-s)\Delta} N^2 (\cdot,s)\| _{L^{\infty}(\Omega)}ds\\[2mm]
\leq&\displaystyle
\mu\int_{t_0}^te^{-r(t-s)} \| N (\cdot,s)\|^2 _{L^{\infty}(\Omega)}ds\\[2mm]\leq
&\mu \varepsilon^2_1\displaystyle\int_{t_0}^t e^{-r(t-s)} e^{-2\alpha (s-t_0)}  ds\\[2mm]
\leq &
\mu \varepsilon^2_1 \displaystyle\int_{t_0}^te^{-(r-\alpha)(t-s)}ds\cdot  e^{-\alpha (t-t_0)}
\\[2mm]
\leq &
 \displaystyle\frac {\mu \varepsilon_1^2} {r-\alpha}e^{-\alpha(t-t_0)}\\
 \leq &\displaystyle\frac {\varepsilon_1} 8
 e^{-\alpha(t-t_0)}
\end{array}
\end{align*}
 due to \eqref{4.2} and $\alpha< r$.  Similarly by \eqref{4.7}, we have
 \begin{align*}
\begin{array}{rl}
J_3\leq
& c_1 \displaystyle\sup_{t\geq t_0}\| u(\cdot,t)\|_{L^\infty(\Omega)} \displaystyle\int_{t_0}^t(1+(t-s)^{-\frac12})e^{-(\lambda_1+r)(t-s)} \|N(\cdot,s)\|_{L^\infty(\Omega)}ds \\[2mm]
\leq &
c_1 \displaystyle\sup_{t\geq t_0}\| u(\cdot,t)\|_{L^\infty(\Omega)}\varepsilon_1 \displaystyle\int_{t_0}^t(1+(t-s)^{-\frac12})e^{-(\lambda_1+r)(t-s)} e^{-\alpha (s-t_0)}
ds \\[2mm]
\leq &
 c_1 \displaystyle\sup_{t\geq t_0}\| u(\cdot,t)\|_{L^\infty(\Omega)}\varepsilon_1 \int_0^\infty(1+\sigma^{-\frac12})e^{-(\lambda_1-\alpha) \sigma}d\sigma \cdot e^{-\alpha (t-t_0)}\\
\leq &
 \displaystyle\frac {\varepsilon_1} 8e^{-\alpha(t-t_0)}.
\end{array}
\end{align*}

As for the term $J_4$, we recall \eqref{4.1}, \eqref{4.8} and apply Lemma 2.1(iii) to get
\begin{align*}
\begin{array}{rl}
J_4\leq
&\chi c_3 \displaystyle\int_{t_0}^t(1+(t-s)^{-\frac23})e^{-(\lambda_1+r)(t-s)} \|(n\nabla w)(\cdot,s)\|_{L^6(\Omega)}ds \\[2mm]
\leq &
\chi c_3\varepsilon_2 \displaystyle\int_{t_0}^t(1+(t-s)^{-\frac23})e^{-(\lambda_1+r)(t-s)}(\frac r \mu+ \varepsilon_1 e^{-\alpha (s-t_0)})  e^{-\alpha (s-t_0)}
ds \\[2mm]
\leq &
\chi c_3\varepsilon_2 (\displaystyle\frac r \mu+\varepsilon_1)\displaystyle\int_{0}^\infty (1+\sigma^{-\frac23})e^{-(\lambda_1+r-\alpha)\sigma}d \sigma\cdot e^{-\alpha (t-t_0)}\\
\leq &
 \displaystyle\frac {\varepsilon_1} 8e^{-\alpha(t-t_0)}+\chi c_3 \displaystyle\frac r \mu I \varepsilon_2 e^{-\alpha (t-t_0)} \\
 \leq &\displaystyle\frac {\varepsilon_1} 4
 e^{-\alpha(t-t_0)}
\end{array}
\end{align*}
provided that  \begin{align}\label{4.21}
\chi c_3\varepsilon_2 \displaystyle\frac r \mu I  < \displaystyle\frac {\varepsilon_1} 8. \end{align}
Therefore, letting $ \varepsilon_2 =4 c_1|\Omega|^{\frac 16}I \varepsilon_1 $ in Lemma 4.1, %and $\tilde{\theta}^*:=$.
 \eqref{4.21} can be warranted  %is then valid
 provided $\mu> 32 \chi c_3c_1|\Omega|^{\frac 16}I^2 r$,  %\tilde{\theta}^*$
 %warrants that
 and thereby
$$
\| N(\cdot,t)\|_{L^\infty(\Omega)}\leq \frac {5\varepsilon_1} 8 e^{-\alpha (t-t_0)} \quad \hbox{for all }\, t\in[t_0,T).
$$
This along with \eqref{4.11} readily shows that $T$ cannot be finite. In combination with the decay property (4.3), a straightford interpolation argument can be employed to prove \eqref{4.5}.
\vskip2mm

%Collecting the previous lemmata directly leads to
% On the basis of the  previous lemmata, we can % Now we are in the position to %prove Theorem 1.1

\it {Proof  of Theorem 1.2.} ~\rm According to \eqref{4.13a} and  $w= -\ln (\frac c{\|c_0\|_{L^\infty(\Omega)}})$, we have
$c(x,t)\leq
\|c_0\|_{L^\infty(\Omega)} e^{-\frac r {2\mu}(t-t_3)}$ for all $t\geq t_3$. On the other hand, if $
\mu_*(\chi, \Omega, r)
:= \max\{\mu_0, 32 \chi c_3c_1|\Omega|^{\frac 16}I^2r\}$, then as an immediate consequence of Theorem 1.1 and
Lemma 4.3, $n(\cdot,t) \rightarrow \frac r\mu$ and  $\frac{|\nabla  c|}c(\cdot,t)\rightarrow0 $
in $L^\infty(\Omega)$ and  $L^6(\Omega)$, respectively,  at an exponential rate  %as $t\rightarrow \infty$
when $\mu >\mu_*(\chi, \Omega, r)$.  Moreover, with the help of the uniform boundedness of $\|\frac{|\nabla c|}c(\cdot,t)\|_{L^\infty(\Omega)}$  with respect to $t>0$, one can show that $\frac{|\nabla c|}c(\cdot,t)\rightarrow0 $ in $L^p(\Omega)$ for any $p>1$ exponentially by the interpolation argument. The proof of this theorem is thus complete.

\subsection{The case $r=0$}
The proof of Theorem 1.3 proceeds on an alternative reasoning. To this end, making use of  %the decay estimates \eqref{2.2} and
the  decay information on $|\nabla w|$ in $L^\infty(\Omega)$ in \eqref{3.61} and the quadratic degradation in the $n-$equation, we first %show the convergence of
%$(u,|\nabla w|)$ with respect to the norm %to $(0,0)$ in $L^\infty(\Omega)$. % as $t\rightarrow \infty$.
%Now, with at hand, %again by utilizing \eqref{2.2},
 turn the decay property of $\|n(\cdot,t)\|_{L^1(\Omega)}$ from  \eqref{2.2} into  an  upper bound estimate of $\|n(\cdot,t)\|_{L^\infty(\Omega)}$. % we can now  obtain the convergence of
%$(u,|\nabla w|)$ with respect to the norm %to $(0,0)$
%in $L^\infty(\Omega)$. % as $t\rightarrow \infty$.
\begin{lemma}\label{lemma45} Let $(n,w,u)$ be the global bounded  solution  of \eqref{2.1} obtained in Theorem 1.1 with $r=0, \mu >0$. Then  one can find constant $C>0$ such that
\begin{equation}\label{4.26}
\| n(\cdot,t)\|_{L^\infty(\Omega)}\leq\frac C{t+1}~~\hbox{for all}~~t>0.
\end{equation}
 \end{lemma}
\it{Proof.}\rm\quad According to the known smoothing properties of  the Neumann heat semigroup $(e^{ \tau\Delta})_{t>0}$ on $\Omega\subset \mathbb{R}^n$ (see \cite{Winkler7}), one can pick $c_1>0$ and $c_2>0$ such that
for all $0<\tau\leq 1$,
\begin{equation}\label{4.27}
\|e^{ \tau\Delta}\varphi\|_{L^\infty(\Omega)}\leq c_1\tau^{-\frac n 2}\|\varphi\|_{L^1(\Omega)}~~\hbox{for all}~~  \varphi\in L^1(\Omega)
\end{equation}
and
\begin{equation}\label{4.28}
\|e^{ \tau\Delta}\nabla\cdot\varphi\|_{L^\infty(\Omega)}\leq c_2\tau^{-\frac 12-\frac n{2p}}\|\varphi\|_{L^p(\Omega)}~~\hbox{for all}~~  \varphi\in C^1(\Omega;\mathbb{R}^n).
\end{equation}
By \eqref{3.63} and \eqref{3.64}, there exists $t_0>3$ such that
\begin{equation}\label{4.29}
24 c_2( \chi\|\nabla w(\cdot,t)\|_{L^3(\Omega)}+   \|u (\cdot,t)\|_{L^2(\Omega)} )\leq 1~~\hbox{ for all} ~~t>t_0-1.
\end{equation}
Now in order to prove the lemma, it is sufficient to derive a bound, independent of $T\in (t_0,\infty)$,  for $M(T)\triangleq \displaystyle \sup_{t_0-1<t<T}\{t\| n(\cdot,t)\|_{L^\infty(\Omega)}\}$.

By the variation-of-constants representation of $n$, we have
\begin{equation}
\begin{array}{ll}
n(\cdot,t)&= e^{\Delta}n(\cdot,t-1)+\displaystyle\chi\!\int_{t-1}^t e^{(t-s)\Delta}\nabla\cdot(n\nabla w)(\cdot,s)ds-
\int_{t-1}^t e^{(t-s)\Delta}(u\cdot\nabla n)(\cdot,s)ds\\
& -
\displaystyle\mu\int_{t-1}^t e^{(t-s)\Delta} n^2(\cdot,s) ds.
\end{array}
\end{equation}
Since $e^{ (t-s)\Delta}$ is nonnegative in $\Omega$ for all $0<s<t$  due to the maximum principle,  it follows from  the
nonnegativity of $n$  that for all $t\in (t_0, T)$
\begin{align*}
&\|n(\cdot,t)\|_{L^\infty(\Omega)}\\
\leq &\| e^{\Delta}n(\cdot,t-1)\|_{L^\infty(\Omega)}+\chi\!\int_{t-1}^t \|e^{(t-s)\Delta}\nabla\cdot(n\nabla w)(\cdot,s)\|
_{L^\infty(\Omega)} ds+
\int_{t-1}^t \|e^{(t-s)\Delta}(u\cdot \nabla n)(\cdot,s)\|
_{L^\infty(\Omega)} ds
\end{align*}
which along with \eqref{4.27}, \eqref{4.28}, \eqref{4.29} and \eqref{2.2} yields
\begin{align*}
&\|n(\cdot,t)\|_{L^\infty(\Omega)}\\
\leq & c_1\|n(\cdot,t-1)\|_{L^1(\Omega)}+c_2\chi\!\int_{t-1}^t (t-s)^{-\frac 56}\|(n\nabla w)(\cdot,s)\|
_{L^3(\Omega)} ds+c_2\!\int_{t-1}^t (t-s)^{-\frac 56}\|(u n)(\cdot,s)\|
_{L^3(\Omega)} ds\\
\leq & \displaystyle\frac {c_1 |\Omega|} {\mu(t-1+\gamma)}+ \frac {6c_2} {t-1} (\chi\displaystyle \max_{t_0-1<s<T}\|\nabla w(\cdot,s)\|
_{L^3 (\Omega)} + \|u (\cdot,t)\|_{L^2(\Omega)} ) \! \cdot M(T)\\
\leq & \displaystyle\frac {c_1 |\Omega|} {\mu(t-1+\gamma)}+ \frac 1 { 4 (t-1)} M(T).
\end{align*}
Hence,
$$
M(T)\leq  \displaystyle\frac {4c_1 |\Omega|} {\mu}+2\displaystyle \sup_{t_0-1<s<t_0}\{s\| n(\cdot,s)\|_{L^\infty(\Omega)}\},
$$
which readily yields \eqref{4.26} since $T>t_0$ is arbitrary, and thus ends the proof. % of this lemma.

In light of Lemma \ref{lemma45}, we can derive  a pointwise estimate $c(x,t)$ from below.
\begin{lemma}\label{lemma46} Let $(n,w,u)$ be the global classical solution  of \eqref{2.1} obtained in Theorem 1.1 with $r=0, \mu >0$. Then  there exists
$\kappa>0$ fulfilling
\begin{align}\label{4.30}
c(x,t)\geq \frac {\displaystyle\inf_{x\in \Omega}c_0(x)}  {(t+1)^{\kappa}}.
\end{align}
\end{lemma}
\it{Proof.}\rm\quad
By the second equation of \eqref{2.1} and Lemma \ref{lemma45}, we can see that
$$w_t\leq \triangle w- |\nabla w|^2 +\frac{c_1}{t+1}-u \cdot\nabla w
$$
with some $c_1>0$ for all $t>0$.  Let $y\in C^1([0,\infty))$ denote the solution of the initial-value problem $y'(t)=\frac{c_1}{t+1}$, $y(0)=\|w_0\|_{L^\infty(\Omega)}$,
then from the  comparison principle, we infer that
\begin{align}\label{4.31}
w(x,t)\leq  \|w_0\|_{L^\infty(\Omega)}+c_1\ln (t+1) ~~\hbox{ for all} ~t>0,
\end{align}
which along with $w= -\ln (\frac c{\|c_0\|_{L^\infty(\Omega)}}) $, yields \eqref{4.30} with $\kappa=c_1$.

Now  utilizing the  decay information on $|\nabla w|$ in $L^\infty(\Omega)$ in \eqref{3.61} again, and thanks to the precise information on the decay of $\| n(\cdot,t)\|_{L^\infty(\Omega)}$ in Lemma \ref{lemma45}, we can obtain the desired  estimate for  $\|n(\cdot,t)\|_{L^\infty(\Omega)}$ from below as well as
the upper estimate for
 $\|\nabla w(\cdot,t)\|_{L^\infty(\Omega)}$.
\begin{lemma}\label{lemma47} Let $(n,w,u)$ be the  solution  of \eqref{2.1} obtained in Theorem 1.1 with $r=0, \mu >0$. Then  one can find
$C_1>0$ and $C_2>0$ fulfilling
\begin{align}\label{4.32}
\| n(\cdot,t)\|_{L^\infty(\Omega)}\geq\frac1{|\Omega|}\| n(\cdot,t)\|_{L^1(\Omega)} \geq \frac{C_1}{t+1}
~~~\hbox{ for all}~~t>0
\end{align}
as well as
\begin{align}\label{4.33}
\|\nabla w(\cdot,t)\|_{L^\infty(\Omega)}\leq\frac{C_2}{t+1}~~\hbox{for  all}~~t>0.
\end{align}
\end{lemma}

\it{Proof.}\rm\quad  We first adapt the  method in  Lemma \ref{lemma45} to derive the precise decay rate of $ \|\nabla w(\cdot,t)\|_{L^\infty (\Omega)} $.
By \eqref{3.63} and \eqref{3.64}, one can choose some  $t_0>2$ such that
\begin{equation}\label{4.34}
4 c_1\int^{\infty}_0 (1+\sigma ^{-\frac 12}) e^{-\lambda_1 \sigma} d\sigma (\|\nabla w(\cdot,t)\|_{L^\infty(\Omega)}+ \|u(\cdot,t)\|_{L^\infty(\Omega)})\leq 1~~\hbox{ for all} ~~t>\frac{t_0}2,
\end{equation}
and then let $M(T)\triangleq \displaystyle \sup_{\frac{t_0}2<s<T}\{s\| \nabla w(\cdot,s)\|_{L^\infty(\Omega)}\}$ for all $T>t_0$.

By the variation-of-constants representation of $w$, we have
$$
 w(\cdot,t)= e^{\frac t2 \Delta}w(\cdot,\frac t2)-\!\int_{\frac t2}^t e^{(t-s)\Delta}|\nabla w|^2(\cdot,s)ds+
\int_{\frac t2}^t e^{(t-s)\Delta}( n-u\cdot\nabla w) (\cdot,s)ds
$$
 for all $t_0<t<T$. We then show that
\begin{align*}
&\|\nabla w(\cdot,t)\|_{L^\infty (\Omega)}
\\
\leq & \|\nabla e^{\frac t2 \Delta}w(\cdot,\frac t2)\|_{L^\infty (\Omega)}
+\!\int_{\frac t2}^t\|\nabla e^{(t-s)\Delta}|\nabla w|^2\|_{L^\infty (\Omega)}+
\int_{\frac t2}^t \|\nabla e^{(t-s)\Delta} n\|_{L^\infty (\Omega)}
\\
& +
\int_{\frac t2}^t \|\nabla e^{(t-s)\Delta} (u\cdot\nabla w)\|_{L^\infty (\Omega)}
\\
\leq & c_1(1+t^{-\frac12})e^{-\frac{\lambda_1t}2} \|w(\cdot,\frac t2)\|_{L^\infty (\Omega)}+ c_1\int_{\frac t2}^t(1+(t-s)^{-\frac 12})e^{-\lambda_1(t-s)} \|n(\cdot,s)\|_{L^\infty (\Omega)}
\\
& +c_1\int_{\frac t2}^t(1+(t-s)^{-\frac 12})e^{-\lambda_1(t-s)} \|\nabla w(\cdot,s)\|_{L^\infty (\Omega)}
(\|\nabla w(\cdot,s)\|_{L^\infty (\Omega)}+\|u(\cdot,s)\|_{L^\infty (\Omega)})
\\
\leq &c_1(1+t^{-\frac12})e^{-\frac{\lambda_1t}2}(\|w_0\|_{L^\infty(\Omega)}+c_2\ln (t+1))+ \frac {2c_1 c_2}{t}
\int^{\infty}_0 (1+\sigma ^{-\frac 12}) e^{-\lambda_1 \sigma}d\sigma
\\
 &+ \frac {2c_1}{t}\int^{\infty}_0 (1+\sigma ^{-\frac 12}) e^{-\lambda_1 \sigma}d\sigma \displaystyle\sup_{t\geq \frac {t_0}2} (\|\nabla w(\cdot,t)\|_{L^\infty(\Omega)}+\|u(\cdot,t)\|_{L^\infty(\Omega)})\cdot M(T)
 \\
 \leq &c_1(1+t^{-\frac12})e^{-\frac{\lambda_1t}2}(\|w_0\|_{L^\infty(\Omega)}+c_2 \ln (t+1))+ \frac {2c_1 c_2}{t}
\int^{\infty}_0 (1+\sigma ^{-\frac 12}) e^{-\lambda_1 \sigma}d\sigma
\\
 &+ \frac {1}{2t} M(T)
\end{align*}
by using Lemma \ref{2.1}(i), \eqref{4.31}, \eqref{4.26} and \eqref{4.34}.
This along with the definition of $M(T)$ yields
$$
M(T)\leq  2\displaystyle \sup_{\frac{t_0}2<s<t_0}\{s\| \nabla w(\cdot,s)\|_{L^\infty(\Omega)}\}
+ c_3
$$ with some constant $c_3>0$ as $\displaystyle\lim_{t\rightarrow \infty} t\ln(t+1) e^{-\lambda_1 t}=0$.
Hence, upon the definition of $M(T)$, we arrive at \eqref{4.33} with an evident choice of $C_2$.

Continuing with the proof, we claim that there exists  $c_4>0$ such that
\begin{align}\label{4.35}
\| n(\cdot,t)\|_{L^\infty(\Omega)}\geq\frac1{|\Omega|}\| n(\cdot,t)\|_{L^1(\Omega)} \geq \frac{c_4}{t+1}
~~~\hbox{ for all}~~t>0.
\end{align}
Indeed, from the $n-$equation of \eqref{2.1} with $r=0$ and  Young's inequality, it follows that
\begin{align*}%\label{4.35}
\displaystyle\frac{d}{dt}\int_\Omega \ln n
&= \displaystyle \int_\Omega  \frac{ |\nabla n|^2}{n^2} +\chi\displaystyle\int_\Omega \frac1 n \nabla\cdot (n\nabla w) - \mu
 \displaystyle\displaystyle\int_\Omega n\nonumber\\
&\geq -\frac{\chi^2}4 \displaystyle \int_\Omega  |\nabla w|^2 - \mu
 \displaystyle\displaystyle\int_\Omega n.
\end{align*}
Inserting \eqref{2.2} and  \eqref{4.33} into the above inequality yields
\begin{align*}
\displaystyle\frac{d}{dt}\int_\Omega \ln n
\geq -\frac{\chi^2}4  \frac{C_2^2|\Omega|}{(t+1)^2} - \displaystyle\frac {|\Omega|} {t+\gamma}
\end{align*}
and thus
\begin{align}\label{4.36}
\displaystyle\int_\Omega \ln n(\cdot,t)
\geq -|\Omega|\ln (t+\gamma)-c_5 ~~~\hbox{ for all}~~t>1
\end{align}
with some $c_5>0$.
On the other hand, by the Jensen inequality, we have
\begin{align*}
|\Omega|\ln (\displaystyle\int_\Omega  n(\cdot,t))-|\Omega|\ln |\Omega|=|\Omega| \ln \{ \frac 1{|\Omega|}\int_\Omega  n(\cdot,t))\}
\geq
\displaystyle\int_\Omega \ln n(\cdot,t).
\end{align*}
This inequality  together with \eqref{4.36} readily leads to \eqref{4.32}. %and thus completes the proof of this lemma.

With the above lemmas at hand, we can now complete the proof of Theorem 1.3.

\it {Proof  of Theorem 1.3.}\rm~~By $w= -\ln (\frac c{\|c_0\|_{L^\infty(\Omega)}}) $, Lemma \ref{lemma45} and Lemma \ref{lemma47}, one can see that $(n,\frac{|\nabla c|}c)\longrightarrow (0,0)
$
in $L^\infty(\Omega)$   algebraically   as $t\rightarrow \infty$. Hence it suffices to show  the decay property of $c(x,t)$. In view of the $w$-equation in \eqref{2.1},
 \eqref{4.35}, \eqref{4.33} and  $\nabla\cdot u=0$,  we can pick  $c_1>0$, $c_2> 0$ and $c_3>0$ such that
\begin{align*}%\label{4.35}
\displaystyle\frac{d}{dt}\int_\Omega w
&=\displaystyle\int_\Omega n -\displaystyle \int_\Omega  |\nabla w|^2- \int_\Omega  u \cdot\nabla w \nonumber\\
&\geq  \frac{c_1|\Omega|}{t+1} - \frac{c_2|\Omega|}{(t+1)^2},~
\end{align*}
 and hence
\begin{align*}
\int_\Omega w(\cdot,t)\geq c_1|\Omega|\ln (t+1)-c_3,
\end{align*}
which entails that for any $t>0$ there exists $x_0(t)\in \Omega$ such that
$$
w(x_0(t),t)\geq c_1\ln (t+1)-\frac {c_3}{|\Omega|}.
$$
Since for  each $\varphi\in W^{1,p}(\Omega)$ with $p>2$, there exists $c_4>0$ such that   $$|\varphi(x)-\varphi(y)|\leq c_4 |x-y|^{1-\frac 2p} \|\nabla \varphi\|_{L^p(\Omega)} ~~\hbox{ for all}~x,y\in \Omega,
$$
we therefore obtain from  \eqref{4.33}  that
\begin{align}
w(x,t) &\geq
w(x_0(t),t)-|x-x_0(t)|\|\nabla w(\cdot,t)\|_{L^\infty(\Omega)}\\
&\geq c_1\ln (t+1)-\frac {c_3}{|\Omega|}-c_4 \hbox{diam}(\Omega),\nonumber
\end{align}
and thereby $$c(x,t)\leq \displaystyle \frac {c_5}{(t+1)^{c_1}}
~~\hbox{ for }~x\in \Omega, t>0
 $$ with some $c_5>0$. This together with \eqref{4.30} shows that $c(x,t)$ actually converges to 0 in $L^\infty(\Omega)$  algebraically  as $t\rightarrow \infty$, and thus ends the proof of Theorem 1.3.

\section*{Acknowledgment}
  Peter Y. H. Pang  was partially supported by the NUS AcRF grant
(R-146-000-249-114). The work of Yifu Wang was supported  by the NNSFC grant  (No.~12071030)  and funded by the Beijing Key Laboratory on MCAACI.
 Jingxue Yin acknowledges support of  NNSFC Grant (No.~11771156), Guangdong Basic and Applied  Basic Research Foundation Grant (No.~2020B1515310013) and NSF of Guangzhou Grant (No.~201804010391).

\end{document}